\documentclass{LMCS}

\usepackage[T1]{fontenc}
\usepackage{amsmath,amssymb,microtype,array,mathrsfs,amsthm,enumerate,mathdefs,graphicx,relsize,hyperref}

\makeatletter
\newcommand\m@thsm@ller[2]{\mbox{\relsize{-1}$\m@th#1#2$}}
\let\smaller\undefined
\DeclareRobustCommand\smaller[1]{\relax\ifmmode{\mathpalette\m@thsm@ller{#1}}\else{\relsize{-1}#1}\fi}
\makeatother

\newenvironment{enuma}{\begin{enumerate}[\upshape(a)]}{\end{enumerate}}
\newenvironment{enumi}{\begin{enumerate}[\upshape(i)]}{\end{enumerate}}
\newenvironment{enumr}{\begin{enumerate}[\upshape(1)]}{\end{enumerate}}

\newcommand*{\dom}{\qopname\relax o{dom}}
\newcommand*{\MSO}{\smaller{\mathrm{MSO}}}
\newcommand*{\GSO}{\smaller{\mathrm{GSO}}}
\newcommand*{\CMSO}{\smaller{\mathrm{CMSO}}}
\newcommand*{\CGSO}{\smaller{\mathrm{CGSO}}}
\newcommand*{\edg}{\mathrm{edg}}
\newcommand*{\WD}{\qopname\relax o{wd}}
\newcommand*{\twd}{\qopname\relax o{twd}}
\newcommand*{\pwd}{\qopname\relax o{pwd}}

\newcommand*{\MTh}{\mathrm{MTh}}
\newcommand*{\qr}{\mathrm{qr}}
\newcommand*{\ar}{\mathrm{ar}}
\newcommand*{\lex}{\mathrm{lex}}
\newcommand*{\ord}{\mathrm{ord}}
\newcommand*{\Min}{\mathrm{Min}}
\newcommand*{\Lf}{\mathrm{Lf}}
\newcommand*{\GF}{\mathrm{Gf}}
\newcommand*{\IN}{\mathrm{in}}
\newcommand\STR{\mathbb{STR}}

\newcommand\TREE{\mathbb{TREE}}
\newcommand*{\emptyseq}{\langle\rangle}
\DeclareRobustCommand*{\Aboveseg}{\mathord\Uparrow}
\newcommand*{\?}{\kern .08em}

\newcommand\nmodels{\not\models}
\newcommand\nsqsubseteq{\not\sqsubseteq}
\newcommand\Sqsubset{\vartriangleleft}

\newcommand\upqed{\vskip-\baselineskip\vskip-\belowdisplayskip}

\def\doi{6 (2:2) 2010}
\lmcsheading%
{\doi}
{1--28}
{}
{}
{Jun.~25, 2008}
{Jun.~22, 2010}
{}   

\begin{document}
\title[Monadic Second-Order Transduction Hierarchy]{On the Monadic Second-Order Transduction Hierarchy}

\author[A.~Blumensath]{Achim Blumensath\rsuper a}
\address{{\lsuper a}TU Darmstadt, Germany}
\email{blumensath@mathematik.tu-darmstadt.de}

\author[B.~Courcelle]{Bruno Courcelle\rsuper b}
\address{{\lsuper b}Institut Universitaire de France and Bordeaux University, LaBRI, France}
\email{courcell@labri.fr}
\thanks{{\lsuper b}Supported by the GRAAL project of the `Agence Nationale pour la Recherche'.}


\begin{abstract}
We compare classes of finite relational structures via monadic
second-order transductions. More precisely, we study the preorder
where we set $C \sqsubseteq K$ if, and only if, there exists a
transduction $\tau$ such that $C\subseteq\tau(K)$.  If we only
consider classes of \emph{incidence structures} we can completely
describe the resulting hierarchy.  It is linear of order type $\omega
+ 3$. Each level can be characterised in terms of a suitable variant
of tree-width.  Canonical representatives of the various levels
are\kern0.08em: the class of all trees of height~$n$, for each $n \in
\bbN$, of all paths, of all trees, and of all grids.
\end{abstract}

\maketitle

\section{Introduction}

\emph{Monadic second-order logic} ($\MSO$) is one of the most expressive logics
for which the theories of many interesting classes of structures
are still decidable.
In particular, the infinite binary tree and many linear orders have
a decidable $\MSO$-theory \cite{Rabin69,Shelah75}
and the same holds for many classes of (finite or infinite)
structures with bounded \emph{tree-width}
\cite{BlumensathColcombetLoeding07,RobertsonSeymour86}.
Furthermore, for every fixed $\MSO$-sentence~$\varphi$ and
every class~$\calC$ of finite structures with bounded tree-width,
there is a linear-time algorithm that checks whether a given
structure from~$\calC$ satisfies~$\varphi$ \cite{Bodlaender96,FlumGrohe06}.
Examples of monadic second-order expressible graph properties are $k$-colourability,
various types of connectivity, and planarity
(via Kuratowski's well-known characterisation by forbidden configurations).

A variant of monadic second-order logic called
\emph{guarded second-order logic} ($\GSO$)
allows quantification not only over sets of elements
but also over sets of edges (i.e., tuples from the relations)~\cite{GraedelHirschOtto02}.
The above mentioned linear-time algorithms can be
adapted to this logic.
There are tight links between guarded second-order logic and tree-width\?:
every class of (finite or infinite) relational structures with a decidable
$\GSO$-theory has bounded tree-width. This gives a sort of converse to the above mentioned
decidability results \cite{Seese91,Courcelle95}. The proof of this result uses a deep
theorem of graph minor theory by Robertson and Seymour\?:
a set of graphs has bounded tree-width if and only if it
excludes some planar graph as a minor~\cite{RobertsonSeymour86}.

To compare the $\MSO$-theories or $\GSO$-theories of two classes of structures
we can use \emph{monadic second-order transductions,}
a certain kind of interpretations suitable both, for monadic second-order logic
and, using a detour via incidence structures, also for guarded second-order logic
\cite{BlumensathColcombetLoeding07,Courcelle91,Courcelle95,Courcelle03}.

In the present article we classify classes
of finite structures according to their `combinatorial complexity'.
(Note that we do not consider decidability issues.)
We will consider two ways to measure the complexity of such classes.
On the one hand, we can use their tree-width and its variants.
On the other hand, we can compare them via transductions.
As it turns out, these two approaches are equivalent and
they give rise to the same hierarchy.
This indicates the robustness of our definitions and their intrinsic interest.
Other possible hierarchies, based on different logics, will be briefly mentioned
in Section~\ref{Sect: conclusion}.

Let us give more details.
An $\MSO$-transduction is a transformation
of relational structures specified by monadic second-order formulae.
As graphs can be represented by relational structures,
we can use $\MSO$-transductions as transformations between graphs.
An $\MSO$-transduction is a generalisation of the following kind
of operations (see Definition~\ref{Def: transduction})\?:

(i) the definition of a relational structure ``inside'' another one
  (in model theory this is called an \emph{interpretation})\?;

(ii) the replacement of a structure~$\frakA$ by the union of a fixed
  number of disjoint copies of~$\frakA$, augmented with appropriate relations
  between the copies\?;

(iii) the expansion of a given structure~$\frakA$ by a fixed number of unary predicates,
  called \emph{parameters.} Usually, these predicates are arbitrary subsets
  of the domain, but we also may have a formula imposing restrictions on them.

Because of the possibility to use parameters, a transduction~$\tau$ is
a many-valued map in general. (We may also think of it as non-deterministic.)
Each relational structure~$\frakA$ has several images
$\tau(\frakA,P_1,\dots,P_n)$ depending on the choice of the parameters
$P_1,\dots,P_n \subseteq A$.
If $\frakB = \tau(\frakA,P_1,\dots,P_n)$ we can consider the tuple~$P_1,\dots,P_n$ as an
\emph{encoding} of~$\frakB$ in~$\frakA$. The transduction~$\tau$ is the corresponding
\emph{decoding function.}

Each transduction~$\tau$ extends in a canonical way to a transformation
between classes of structures.
If $\calC$~and~$\calK$ are classes of relational structures with
$\calC \subseteq \tau(\calK)$, we can think of~$\tau$ as a way of
encoding the structures in~$\calC$ by elements of~$\calK$.
For instance, every finite graph can be encoded in a sufficiently large finite square grid
(by a fixed transduction~$\tau$). Every finite tree of height at most~$n$ (for fixed~$n$)
can be encoded in a sufficiently long finite path.
But it is not the case that all finite trees can be encoded by paths (by a single transduction).

The purpose of this article is to classify classes of finite relational structures
according to their \emph{encoding powers.}
We will compare classes $\calC$~and~$\calK$ of structures by the following preorder\?:
\begin{align*}
  \calC \sqsubseteq \calK \quad\defiff\quad
  \calC \subseteq \tau(\calK) \quad\text{for some $\MSO$-transduction } \tau\,.
\end{align*}
We attack the problem of determining the structure of this preorder.
Since, at the moment, a complete description of this hierarchy seems to be out of reach,
we concentrate on a variant where we replace monadic second-order logic by
\emph{guarded second-order logic.}
In this case, the corresponding hierarchy can be described completely.
To obtain a corresponding notion of transduction we cannot simply change
the definition of an $\MSO$-transduction to use $\GSO$-formulae since
the resulting notion of transduction would not yield a reduction between $\GSO$-theories,
and it even would not be closed under composition.
Instead, we will take a detour by combining ordinary $\MSO$-transductions
with a well-known translation between $\GSO$ and $\MSO$.

This translation is based on \emph{incidence structures.}
Let us first describe this notion for undirected graphs
where it is very natural.
There are two canonical ways to encode a graph~$\frakG$ by a relational structure.
We can use its \emph{adjacency representation} which is a structure
$\langle V,\mathrm{edg}\rangle$ where the domain~$V$ consists of all vertices
of~$\frakG$ and $\mathrm{edg}$ is a binary relation containing all pairs of adjacent
vertices. But we also can use the \emph{incidence representation} of~$\frakG$.
This is the structure $\langle V \cup E, \mathrm{in}\rangle$ where the domain
$V \cup E$ contains both, the vertices and the edges of~$\frakG$, and
$\mathrm{in}$ is the incidence relation between vertices and edges.
In a similar way, we can associate with every relational structure~$\frakA$
its incidence structure~$\frakA_\IN$ (see Definition~\ref{Def: incidence structure})
where the domain also contains elements for all tuples in some relation of~$\frakA$.

It is shown in \cite{GraedelHirschOtto02} that every $\GSO$-formula~$\varphi$
talking about some structure~$\frakA$ can be translated into an $\MSO$-formula
talking about the incidence structure~$\frakA_\IN$, and vice versa.
Hence, we can use incidence structures to obtain an analogue~$\sqsubseteq_\IN$
of the preorder~$\sqsubseteq$ suitable for guarded second-order logic. We set
\begin{align*}
  \calC \sqsubseteq_\IN \calK \quad\defiff\quad
  \calC_\IN \subseteq \tau(\calK_\IN) \quad\text{for some $\MSO$-transduction } \tau\,,
\end{align*}
where $\calC_\IN := \set{ \frakA_\IN }{ \frakA \in \calC }$.
The main result of the present article is a complete characterisation of the resulting
hierarchy for classes of finite structures.
We show that the preorder~$\sqsubseteq_\IN$ is linear
of order type $\omega + 3$.
It turns out that every class of finite structures is equivalent
to one of the following classes, listed in increasing order of generality\?:
\begin{itemize}
\item trees of height at most~$n$, for each $n \in \bbN$\?;
\item paths\?;
\item arbitrary trees (equivalently, binary trees)\?;
\item (square) grids.
\end{itemize}
Each of these levels can be characterised in terms of tree decompositions.
Hence, we also obtain a corresponding hierarchy of complexity measures
on structures that are compatible with $\MSO$-transductions transforming
incidence structures.

The upper levels of the hierarchy can be determined easily using techniques
from graph minor theory developed by Robertson and Seymour, such as
the notions of a minor and a tree decomposition.
In particular, we employ two results characterising bounded tree-width and bounded
path-width in terms of excluded minors \cite{RobertsonSeymour83,RobertsonSeymour86}.

For the lower levels, which consist of classes of bounded path-width, the characterisation
is more complicated and requires new results relating tree decompositions
and monadic second-order logic.

In Sections \ref{Sect: preliminaries}~and~\ref{Sect: MSO}
we give basic definitions.
Section~\ref{Sect: tree decompositions} collects some known results from
graph minor theory. We also introduce a new variant of tree-width
and prove some of its basic properties.
In Section~\ref{Sect: tree decompositions and transductions}
we expound the connections between tree-width
and monadic second-order transductions.
In Section~\ref{Sect: hierarchy} we introduce the transduction hierarchy
and we state our main theorem.
Its proof is contained in
Sections \ref{Sect: hierarchy2}~and~\ref{Sect: hierarchy3}.
In the first one, we prove that the hierarchy is strict
while, in the second one, we show that it covers every class.
The final Section~\ref{Sect: conclusion} contains some extension
of our results to other logics and some open problems in this direction.

\section{Preliminaries}
\label{Sect: preliminaries}

Let us fix our notation.
We set $[n] := \{0,\dots,n-1\}$
and we write $\PSet(X)$ for the power set of a set~$X$.
We denote tuples~$\bar a$ with a bar. The components of~$\bar a$
will be $a_0,\dots,a_{n-1}$ where the length~$n$ will usually be implicit.
We sometimes identify a tuple~$\bar a$ with the set of its components.
For instance, we write $c \in \bar a$ to express that $c = a_i$, for some~$i$.

In this article all graphs, trees, and relational structures are finite.
We will not repeat this finiteness assumption.
A relational structure~$\frakA$ is of the form
$\langle A,R_0^\frakA,\dots,R_{m-1}^\frakA\rangle$
with domain~$A$ and relations~$R_i^\frakA$.
The \emph{signature} of such a structure is the set $\Sigma = \{R_0,\dots,R_{m-1}\}$
of relation symbols.
In some proofs we will also use signatures with constant symbols denoting elements of the domain.
We write $\ar(R)$ for the arity of a relation~$R$.
For a signature~$\Sigma$, we denote by $\STR[\Sigma]$ the class of all $\Sigma$-structures.
We write $\frakA \oplus \frakB$ for the \emph{disjoint union} of the structures $\frakA$~and~$\frakB$.

We mainly consider \emph{incidence structures.} These are representations
of structures~$\frakA$ where we have added new elements to the domain, one for each
tuple in the relations of~$\frakA$.
\begin{defi}\label{Def: incidence structure}
Let $\frakA = \langle A,R^\frakA_0,\dots,R^\frakA_{m-1}\rangle$ be a structure
and let $r$~be the maximal arity of a relation~$R_i$.
The \emph{incidence structure} of~$\frakA$ is the structure
\begin{align*}
  \frakA_\IN :=
    \langle A \cupdot E, P_{R_0},\dots,P_{R_{m-1}}, \IN_0,\dots,\IN_{r-1}\rangle\,,
\end{align*}
where we extend the domain~$A$ by
\begin{align*}
  E := R^\frakA_0 \cup\dots\cup R^\frakA_{m-1}\,,
\end{align*}
and the relations are
\begin{align*}
  P_{R_i} &:= \set{ \bar c \in E }{ \bar c \in R_i^\frakA }\,, \\
  \IN_i &:= \set{ (a,\bar c) \in A \times E }{ \abs{\bar c} > i \text{ and } a = c_i }\,.
\end{align*}
The class of all incidence structures is $\STR_\IN[\Sigma] := \set{ \frakA_\IN }{ \frakA \in \STR[\Sigma] }$.
\end{defi}
\begin{rem}
Note that incidence structures are binary (i.e., their relations have arity at most~$2$).
Hence, they can be regarded as bipartite labelled directed graphs.
\end{rem}

One important property of incidence structures is the fact that they are
\emph{sparse,} i.e., their relations contain few tuples.
\begin{defi}
Let $k \in \bbN$. A structure\/ $\frakA = \langle A,\bar R\rangle$ is \emph{$k$-sparse}%
\footnote{\cite{Courcelle03} introduced two notions of sparsity for hypergraphs\?:
\emph{$k$-sparse hypergraphs} and \emph{uniformly $k$-sparse hypergraphs.}
What we call $k$-sparse is a slight modification of the uniform version of \cite{Courcelle03}.}
if, for every subset $X \subseteq A$ and all relations~$R_i$, we have
\begin{align*}
  \bigabs{R_i \cap X^{\ar(R_i)}} \leq k\cdot\abs{X}\,.
\end{align*}
\end{defi}
\begin{lem}
Every incidence structure is $1$-sparse.
\end{lem}

Let us fix our notation regarding trees and graphs.
\begin{defi}
A \emph{directed graph} is a pair $\langle V,\edg\rangle$
where $V$~is the set of \emph{vertices} and $\edg \subseteq V \times V$
is the \emph{edge relation.}
Thus, graphs are by definition simple (without parallel edges).
An \emph{undirected graph} is a graph where the edge relation $\edg$
is symmetric. When speaking of a \emph{graph} we will always mean an undirected one.

We regard a \emph{coloured graph} as a relational structure
$\langle V, E_0,\dots,E_k,P_0,\dots,P_m\rangle$
with binary relations~$E_i$ and unary relations~$P_i$ that encode
the colours of, respectively, the edges and the vertices.
We allow graphs whose edges and vertices have several colours.
\end{defi}

Trees play a major role in this article.
Intuitively, a tree is a directed graph~$\frakT$
with a unique vertex~$r$ of indegree~$0$,
called the \emph{root} of~$\frakT$,
such that every vertex is reachable from~$r$
by a unique directed path.
The actual definition we will use
is slightly more concrete.
It is based on the usual encoding of the vertices
of a tree by finite sequences describing the path from
the root to the given vertex.
In fact, we introduce two notions of a tree\?:
\emph{order-trees} and \emph{successor-trees.}
The latter use the usual edge relation, while
the former are equipped with the tree-order instead.
\begin{defi}
Let $D$~be a set.

(a) For an ordinal~$\alpha$, we denote by $D^{<\alpha}$
the set of all sequences of elements of~$D$ of length less than~$\alpha$.
The \emph{prefix relation} on~$D^{<\omega}$ is defined by
\begin{align*}
  x \preceq y \quad\defiff\quad y = xz\,, \quad\text{for some } z \in D^{<\omega}.
\end{align*}
The infimum of $x$~and~$y$ with respect to~$\preceq$, i.e., their longest common prefix,
is denoted by $x \sqcap y$.

(b) A finite prefix closed subset $T \subseteq D^{<\omega}$ is called a \emph{tree domain.}
Following our intuition that a vertex is represented by the path leading to it,
we call the empty sequence~$\emptyseq$ the \emph{root} of~$T$
and the maximal elements of~$T$ its \emph{leaves.}
The domain of the \emph{complete $m$-ary tree of height~$n$} is $m^{<n}$.
Hence, $m^{<0}$~is the empty tree, $m^{<1}$~the one consisting only of the root,
and $m^{<2}$ consists of a root and $m$~leaves.

Given a tree domain~$T$ we can define the \emph{successor relation}
$\edg$ on~$T$ by setting
\begin{align*}
  \langle x,y\rangle \in \edg \quad\defiff\quad
  y = xd \text{ for some } d \in D\,.
\end{align*}
In this case we call $y$~a \emph{successor} of~$x$
and $x$~the \emph{predecessor} of~$y$.
A structure of the form $\langle T,\edg\rangle$ (and every structure
isomorphic to it) is called a \emph{successor-tree.}

Sometimes it is convenient to replace the successor relation~$\edg$
by the tree order~$\preceq$. Structures of the form $\langle T,{\preceq}\rangle$
are called \emph{order-trees.}
A~\emph{coloured tree} is the expansion of a (order- or successor-) tree
by unary predicates~$\bar P$.
(We do not require these predicates to be pairwise disjoint.
Hence, every vertex may have none, one, or several colours.)
We write $\TREE_m$ for the class of all order-trees
$\langle T,{\preceq},P_0,\dots,P_{m-1}\rangle$ with $m$~colours.
The set of \emph{leaves} of a tree~$\frakT$ is denoted by $\Lf(\frakT)$.

(c) Let $\frakT = \langle T,{\preceq}\rangle$ be an order-tree.
The \emph{level} of an element $v \in T$ is the number of vertices $u \in T$
with $u \prec v$. We denote it by~$\abs{v}$.
The \emph{height} of~$\frakT$ is the least ordinal~$\alpha$ greater than
the level of every element of~$T$.
Hence, the empty tree has height~$0$ and the tree with a single vertex has height~$1$.
The \emph{out-degree} of~$\frakT$ is the maximal number of successors
of a vertex of~$\frakT$.
For successor-trees we define these notions analogously.

(d) Let $\frakT$~be a tree and $v$~a vertex of~$\frakT$.
The \emph{subtree} of~$\frakT$ \emph{rooted} at~$v$
is the subtree~$\frakT_v$
consisting of all vertices~$u$ with $v \preceq u$, i.e.,
all vertices below~$v$.
\end{defi}

Sometimes it is possible to reduce statements about relational structures to
statements about graphs. One way to do so consists in replacing
a structure by its \emph{Gaifman graph.}
\begin{defi}
The \emph{Gaifman graph} of a structure $\frakA = \langle A,\bar R\rangle$
is the undirected graph
\begin{align*}
  \GF(\frakA) := \langle A,\edg\rangle\,,
\end{align*}
with the same domain~$A$ and with the edge relation
\begin{align*}
  \edg := \set{ (u,v) }
              { u \neq v \text{ and there is some } \bar c \in R^\frakA_i \text{ with } u,v \in \bar c }\,.
\end{align*}
\end{defi}

\section{Monadic second-order logic and transductions}
\label{Sect: MSO}

Monadic second-order logic ($\MSO$) is the extension of first-order logic
by set variables and quantifiers over such variables.
An important variant of $\MSO$ is
guarded second-order logic ($\GSO$) where one can quantify not only over sets of elements
but also over sets of tuples from the relations (see~\cite{GraedelHirschOtto02} for details).
Hence, guarded second-order logic over a given structure~$\frakA$ is equivalent to
monadic second-order logic over its incidence structure~$\frakA_\IN$.
\begin{lem}[\cite{GraedelHirschOtto02}]\leavevmode\\
\itm{(a)} For every $\GSO$-sentence~$\varphi$, we can effectively construct an $\MSO$-sentence~$\psi$ such that
\begin{align*}
  \frakA \models \varphi \quad\iff\quad \frakA_\IN \models \psi\,,
  \quad\text{for all structures } \frakA\,.
\end{align*}

\noindent
\itm{(b)} For every $\MSO$-sentence~$\varphi$, we can effectively construct a $\GSO$-sentence~$\psi$ such that
\begin{align*}
  \frakA_\IN \models \varphi \quad\iff\quad \frakA \models \psi\,,
  \quad\text{for all structures } \frakA\,.
\end{align*}
\end{lem}
Throughout the article we will consistently work with incidence structures,
thereby avoiding the treatment of guarded second-order logic.
In particular, all formulae are tacitly assumed to be $\MSO$-formulae.

Besides $\MSO$ and $\GSO$ we also consider their \emph{counting extensions}
$\CMSO$ and $\CGSO$. These add predicates of the form $\abs{X} \equiv k \pmod m$
to, respectively, $\MSO$ and $\GSO$, where $X$~is a set variable and $k,m$ are numbers.
All of our results for $\GSO$ go through also for $\CGSO$, i.e., for
$\CMSO$-transductions between incidence structures.
In Section~\ref{Sect: conclusion} we will give a partial characterisation
of the hierarchy for $\CMSO$-transductions of graphs (not of incidence graphs).
In this case the availability of counting predicates does make a difference.

To state the composition theorem below it is of advantage to work with a
variant of $\MSO$ without first-order variables. This variant has atomic
formulae of the form $X \subseteq Y$ and $R\bar Z$, for set variables $X,Y,Z_0,Z_1,\dots$,
where a formula of the form~$R\bar Z$ states that there are elements $a_i \in Z_i$
such that the tuple~$\bar a$ is in~$R$.
Note that every general monadic second-order formula with first-order variables
can be brought into this restricted form
by replacing all first-order variables by set variables and adding the condition
that these sets are singletons.

Whenever we speak of $\MSO$ we will have this version in mind. In particular,
the following definition of the rank of a formula is based on this variant.
When writing down concrete formulae, on the other hand, we will allow the use of first-order
variables to improve readability. We regard every such formula as an abbreviation
of a formula of the restricted form.
Similarly, when we use structures with constants we actually regard each constant
as a singleton set.
\begin{defi}
(a) The \emph{rank}~$\qr(\varphi)$ of a formula~$\varphi$
is the nesting depth of quantifiers in~$\varphi$.
Formulae of rank~$0$ are called \emph{quantifier-free.}

(b) The \emph{monadic theory of rank~$m$} of a structure~$\frakA$ is
\begin{align*}
  \MTh_m(\frakA) := \set{ \varphi \in \MSO }{ \frakA \models \varphi,\ \qr(\varphi) \leq m }\,.
\end{align*}
For a tuple~$\bar a$ of elements of~$\frakA$, we also consider the monadic theory
$\MTh_m(\frakA,\bar a)$ of the expansion $\langle\frakA,\bar a\rangle$.
\end{defi}
\begin{rem}
We use the term `rank' instead of the more natural `quantifier rank' since
in Section~\ref{Sect: conclusion} below we will consider $\CMSO$
where the notion of rank has to be adapted for our results
to go through.
\end{rem}

In order to compare the monadic theories of two classes of structures
we employ $\MSO$-transductions.
To simplify the definition we introduce three simple operations
and we obtain $\MSO$-transductions as compositions of these.
\begin{defi}\label{Def: transduction}
(a) Let $k \geq 2$ be a natural number. The operation $\mathrm{copy}_k$ maps a
structure~$\frakA$ to the expansion
\begin{align*}
  \mathrm{copy}_k(\frakA) := \langle \frakA\oplus\dots\oplus\frakA, {\sim}, P_0,\dots,P_{k-1}\rangle
\end{align*}
of the disjoint union of $k$~copies of~$\frakA$ by the following relations.
Denoting the copy of an element $a \in A$ in the $i$-th component of $\frakA \oplus\dots\oplus \frakA$
by the pair $\langle a,i\rangle$, we define
\begin{align*}
  P_i := \set{ \langle a,i\rangle }{ a \in A }
  \qtextq{and}
  \langle a,i\rangle \sim \langle b,j\rangle \ \defiff\ a = b\,.
\end{align*}
For $k = 1$, we set $\mathrm{copy}_1(\frakA) := \frakA$.

(b) For $m \in \bbN$, we define the operation $\mathrm{exp}_m$ that maps
a structure~$\frakA$ to all possible expansions by $m$~unary predicates
$Q_0,\dots,Q_{m-1} \subseteq A$.
Note that this operation is many-valued and that $\mathrm{exp}_0$ is just the
identity.

(c) A \emph{basic $\MSO$-transduction} is a partial operation~$\tau$
on relational structures described by a list
\begin{align*}
  \bigl\langle \chi, \delta(x), \varphi_0(\bar x),\dots,\varphi_{s-1}(\bar x)\bigr\rangle
\end{align*}
of $\MSO$-formulae called the \emph{definition scheme} of~$\tau$.
Given a structure~$\frakA$ that satisfies the formula~$\chi$
the operation~$\tau$ produces the structure
\begin{align*}
  \tau(\frakA) := \langle D, R_0,\dots,R_{s-1}\rangle
\end{align*}
where
\begin{align*}
  D := \set{ a \in A }{ \frakA \models \delta(a) }
  \qtextq{and}
  R_i := \set{ \bar a \in D^{\ar(R_i)} }{ \frakA \models \varphi_i(\bar a) }\,.
\end{align*}
If $\frakA \nmodels \chi$ then $\tau(\frakA)$ remains undefined.

(d) A \emph{$k$-copying $\MSO$-transduction}~$\tau$ is a
(many-valued) operation on relational structures of the form
$\tau_0 \circ \mathrm{copy}_k \circ \mathrm{exp}_m$
where $\tau_0$~is a basic $\MSO$-trans\-duc\-tion.
When the value of~$k$ does not matter, we will simply speak of
a \emph{transduction.}

Note that, due to~$\mathrm{exp}_m$, a structure can be mapped
to several structures by a transduction.
Consequently, we consider $\tau(\frakA)$ as the \emph{set} of possible values
$(\tau_0 \circ \mathrm{copy}_k)(\frakA, \bar P)$ where $\bar P$~ranges
over all $m$-tuples of subsets of~$A$.

For a class~$\calC$, we set
\begin{align*}
  \tau(\calC) := \bigcup {\set{ \tau(\frakA) }{ \frakA \in \calC }}\,.
\end{align*}
\end{defi}

\begin{rem}
(a) The expansion by $m$~unary predicates corresponds, in the terminology
of \cite{Courcelle95,Courcelle03}, to using $m$~\emph{parameters.}
We will use this terminology, for instance,
in the proof of Theorem~\ref{Thm: interpretable implies tree decomposable}.

(b) Note that every basic $\MSO$-transduction is a $1$-copying $\MSO$-transduction
without parameters.
\end{rem}
\begin{exa}\label{exam: gf interpretable}
(a) Let $\Sigma$~be a signature and let $r$~be the maximal arity of a relation in~$\Sigma$.
The operation mapping an incidence structure $\frakA_\IN \in \STR_\IN[\Sigma]$
to the structure $\GF(\frakA)_\IN$ is a $k$-copying $\MSO$-transduction
where $k = r(r-1)/2$, for $r \geq 2$, and $k = 1$, for $r \leq 1$.

(b) For every fixed number~$n \in \bbN$, we describe a
transduction~$\tau$ transforming a path~$\frakP$ of length~$l$ into the
class of all trees of height~$n$ with $l+1$ vertices.

We can encode a tree~$T$ of height~$n$ with $m$~vertices
as a finite word~$w$ of length~$m$ over the alphabet~$[n]$ as follows.
Let $v_0 <_\lex\dots<_\lex v_{m-1}$ be the enumeration of the vertices of~$T$
in lexicographic order, and let $l_i$~be the level of~$v_i$.
We encode~$T$ by the word $w := l_0\dots l_{m-1}$.
A transduction can recover~$T$ from~$w$ as follows.
Each position in~$w$ corresponds to a vertex.
The predecessor of the $i$-th vertex~$v$ is the maximal vertex to the left of~$v$
whose label is less than~$l_i$.
Clearly, this predecessor relation is definable in monadic second-order logic.
\end{exa}

The two most important properties of $\MSO$-transductions
are summarised in the following lemmas.
\begin{lem}\label{Lem: transductions are comorphisms}
Let $\tau$~be a transduction. For every $\MSO$-sentence~$\varphi$, there
exists an $\MSO$-sentence~$\varphi^\tau$ such that, for all structures~$\frakA$,
\begin{align*}
  \frakA \models \varphi^\tau
  \quad\iff\quad
  \frakB \models \varphi \quad\text{for some } \frakB \in \tau(\frakA)\,.
\end{align*}
Furthermore, if $\tau$~is quantifier-free, then the rank of~$\varphi^\tau$
is no larger than that of~$\varphi$.
\end{lem}
\begin{cor}\label{Cor: transductions operate on theories}
For every quantifier-free transduction~$\tau$ and every $m \in \bbN$,
there exists a function~$f_\tau$ on monadic theories of rank~$m$ such that
\begin{align*}
  \MTh_m(\tau(\frakA)) = f_\tau\bigl(\MTh_m(\frakA)\bigr)\,,
  \qquad\text{for all structures } \frakA\,.
\end{align*}
\end{cor}
\begin{lem}[\cite{Courcelle91}]
For all transductions $\sigma,\tau$ there exists a transduction~$\varrho$
such that $\varrho = \sigma \circ \tau$.
\end{lem}

As a further example note that we can use transductions
to translate between order-trees and successor-trees.
\begin{lem}
\itm{(a)}
There exists a transduction~$\tau$ mapping
an order-tree to the corresponding successor-tree.

\itm{(b)}
There exists a transduction~$\sigma$ mapping
a successor-tree to the corresponding order-tree.
\end{lem}

Similarly there are transductions
translating between a structure and its incidence structure.
\begin{lem}
For every signature~$\Sigma$, there exists a transduction~$\tau$
such that $\tau(\frakA_\IN) = \frakA$,
for all $\frakA \in \STR[\Sigma]$.
\end{lem}
The converse statement is a much deeper result
and requires the structure in question to be $k$-sparse
for some fixed~$k$.
\begin{thm}[\cite{Courcelle03,Blumensath10}]
For every signature~$\Sigma$ and all numbers $k \in \bbN$,
there exists an $\MSO$-trans\-duc\-tion~$\tau$
such that $\tau(\frakA) = \frakA_\IN$,
for all $k$-sparse structures $\frakA \in \STR[\Sigma]$.
\end{thm}

We have seen in Lemma~\ref{Lem: transductions are comorphisms}
that transductions relate the monadic theories of two structures.
We also need techniques to relate the monadic theory of a structure
to those of its substructures.
The disjoint union operation can frequently be used for this purpose
(for a proof of the following theorem see, e.g., Theorem 7.11 of~\cite{Libkin04},
or~\cite{Courcelle87}).
\begin{thm}\label{Thm: composition for unions}
Let $\Sigma$~and~$\Gamma$ be relational signatures with constants.
For every $m \in \bbN$, there exists a \textup(computable\textup) binary
operation~$\oplus_m$ on monadic theories of rank~$m$ such that
\begin{align*}
  \MTh_m(\frakA \oplus \frakB) = \MTh_m(\frakA) \oplus_m \MTh_m(\frakB)\,,
\end{align*}
for all $\Sigma$-structures~$\frakA$ and $\Gamma$-structures~$\frakB$.
\end{thm}

Below we will mainly make use of the following corollary.
\begin{lem}\label{Lem: composition for trees}
Let\/ $\frakT$~be an order-tree and $v \in T$ a vertex.
Suppose that\/ $\frakT'$~is the order-tree obtained from\/~$\frakT$ by replacing
the subtree\/~$\frakT_v$ by some tree\/~$\frakS$. Let $\bar c$~be a tuple of vertices of\/~$\frakT$
with $v \npreceq c_i$, for all~$i$. If $\bar a$~are vertices of\/~$\frakT_v$ and $\bar b$~are vertices
of\/~$\frakS$ such that
\begin{align*}
  \MTh_m(\frakT_v,\bar a) = \MTh_m(\frakS,\bar b)
\end{align*}
then it follows that
\begin{align*}
  \MTh_m(\frakT,\bar a\bar c) = \MTh_m(\frakT',\bar b\bar c)\,.
\end{align*}
\end{lem}
\begin{Proof}
Let $\frakC$~be the tree obtained from~$\frakT$ by replacing the
subtree~$\frakT_v$ by a single vertex~$w$. We define the following auxiliary predicates\?:
\begin{align*}
  P := \{w\}\,,\quad
  Q_0 := \set{ x \in C }{ x \prec w }\,,\quad
  Q_1 := T_v\,,
  \qtextq{and}
  Q'_1 := S\,.
\end{align*}
We construct a quantifier-free transduction~$\tau$ such that
\begin{align*}
  \tau\bigl(\langle\frakC,P,Q_0,\bar c\rangle \oplus \langle\frakT_v,Q_1,\bar a\rangle\bigr)
  &= \langle\frakT,\bar a\bar c\rangle\,, \\
  \tau\bigl(\langle\frakC,P,Q_0,\bar c\rangle \oplus \langle\frakS,Q'_1,\bar b\rangle\bigr)
  &= \langle\frakT',\bar b\bar c\rangle\,.
\end{align*}
If $f_\tau$~is the function from Corollary~\ref{Cor: transductions operate on theories}
and $\oplus_m$~the operation from Theorem~\ref{Thm: composition for unions}, it follows that
\begin{align*}
  \MTh_m(\frakT',\bar b\bar c)
  &= f_\tau\bigl(\MTh_m(\frakC,P,Q_0,\bar c) \oplus_m \MTh_m(\frakS,Q'_1,\bar b)\bigr) \\
  &= f_\tau\bigl(\MTh_m(\frakC,P,Q_0,\bar c) \oplus_m \MTh_m(\frakT_v,Q_1,\bar a)\bigr)
   = \MTh_m(\frakT,\bar a\bar c)\,,
\end{align*}
as desired.

Hence, it remains to define~$\tau$.
Let $\{{\preceq},P,Q_0,\bar d\}$ be the signature of
$\langle\frakC,P,Q_0,\bar c\rangle$,
and $\{{\preceq},Q_1,\bar e\}$ the signature
of $\langle\frakT_v,Q_1,\bar a\rangle$ and $\langle\frakS,Q'_1,\bar b\rangle$.
For $\tau$~we can use the basic $\MSO$-transduction consisting of the following formulae\?:
\begin{align*}
  \chi                   &:= \mathrm{true}\,, \\
  \delta(x)              &:= \neg Px\,, \\
  \varphi_{\preceq}(x,y) &:= x \preceq y \lor (Q_0x \land Q_1y)\,.
\end{align*}
\upqed
\end{Proof}

\section{Minors and tree decompositions}   
\label{Sect: tree decompositions}

Some properties of the transduction hierarchy,
which we will introduce in Section~\ref{Sect: hierarchy} below,
can be deduced from results about graph minors.
\begin{defi}
(a) Let $\frakG = \langle V,\edg\rangle$ be an undirected graph and $E \subseteq \edg$
a set of edges. We denote by $E^*$~the reflexive and transitive closure of~$E$.
Note that $E^*$~is an equivalence relation.
The graph $\frakG/E$ is obtained by contracting all edges in~$E$.
Formally, we have
\begin{align*}
  \frakG/E := \langle W,\edg_0\rangle\,,
\end{align*}
where $W := V/E^*$ is the set of equivalence classes
and $\edg_0$~contains an edge between classes $[x]$ and $[y]$
if and only if
$[x] \neq [y]$ and there are representatives $u \in [x]$ and $v \in [y]$ with
$\langle u,v\rangle \in \edg$.

(b) A \emph{minor} of a graph~$\frakG$ is a graph that can be obtained
from~$\frakG$ by first deleting some vertices and edges and then contracting some
of the remaining edges.
For a class~$\calC$ of graphs, we denote by $\Min(\calC)$ the class
of all minors of graphs in~$\calC$.
\end{defi}

One central tool in graph minor theory is the notion of a tree decomposition
and the related complexity measures called tree-width and path-width.
These notions extend in a natural way to relational structures.
\begin{defi}
Let $\frakA = \langle A,\bar R\rangle$ be a structure.

(a) A \emph{tree decomposition} of~$\frakA$ is a family $D = (U_v)_{v \in T}$
of (possibly empty) subsets $U_v \subseteq A$ indexed by a rooted tree~$T$ such that
\begin{itemize}
\item for every element $a \in A$, the set $\set{ v \in T }{ a \in U_v }$
  is nonempty and connected in~$T$\?;
\item for every tuple $\bar c \in R_i$, there is some index $v \in T$ with
  $\bar c \subseteq U_v$.
\end{itemize}
We call the sets~$U_v$ the \emph{components} of the decomposition,
and $T$~is its \emph{underlying tree.}

The \emph{height} of a tree decomposition $D = (U_v)_{v \in T}$ is the height
of~$T$, while its \emph{width} is the number
\begin{align*}
  \WD(D) := \sup_{v \in T} {(\abs{U_v}-1)}\,.
\end{align*}

(b) The \emph{tree-width} $\twd(\frakA)$ of~$\frakA$ is the minimal width
of a tree decomposition of~$\frakA$.

(c) The \emph{path-width} $\pwd(\frakA)$ of~$\frakA$ is the minimal width
of a tree decomposition of~$\frakA$ where the underlying tree is a path.

(d) The \emph{$n$-depth tree-width} $\twd_n(\frakA)$ of~$\frakA$
is the minimal width of a tree decomposition
of~$\frakA$ whose underlying tree has height at most~$n$.

(e) For a class~$\calC$ of structures, we define $\twd(\calC)$
as the supremum of $\twd(\frakA)$, for $\frakA \in \calC$,
and similarly for $\pwd(\calC)$ and $\twd_n(\calC)$.
\end{defi}
\begin{rem}
(a) The $n$-depth tree-width of a graph~$\frakG$
is related to its \emph{tree-depth} $\mathrm{td}(\frakG)$
as introduced by Ne\v set\v ril and Ossona de Mendez
\cite{NesetrilOssona06a,NesetrilOssona06b}.
The tree-depth of a graph~$\frakG$ is the least number~$n$
such that some orientation of~$\frakG$ is a subgraph of some order-tree of height~$n$.
For a graph~$\frakG$, it follows that
\begin{itemize}
\item $\mathrm{td}(\frakG) \leq n$ implies $\twd_n(\frakG) < n$\?;
\item $\twd_n(\frakG) < k$ implies $\mathrm{td}(\frakG) \leq nk$.
\end{itemize}
These facts are easy to establish.
We will not need them in the following.

(b) There are some simple relations between $n$-depth tree-width,
path-width, and tree-width.
For every graph~$\frakG$, we have
\begin{alignat*}{-1}
  \twd(\frakG) &\leq \twd_{n+1}(\frakG) \leq \twd_n(\frakG)\,,
  &&\quad\text{for every } n \in \bbN\,, \\
\prefixtext{and}
  \twd(\frakG) &= \twd_n(\frakG)\,,
  &&\quad\text{for all sufficiently large } n \in \bbN\,. \\
\intertext{Furthermore,}
  \pwd(\frakG) &< n (\twd_n(\frakG)+1)\,,
  &&\quad\text{for every } n \in \bbN\,.
\end{alignat*}
(Let us sketch the proof of the last inequality\?:
let $(U_v)_{v \in T}$ be a tree decomposition of~$\frakG$ of height~$n$
and width $\twd_n(\frakG)$.
As the components of a path decomposition of~$\frakG$ we take
all sets of the form $U_{v_0} \cup\dots\cup U_{v_k}$,
where $v_0\dots v_k$ is a path from the root~$v_0$ to some leaf~$v_k$ of~$T$.)
\end{rem}

The next lemma shows that
most questions regarding tree decompositions of a structure can
be reduced to the corresponding questions about its Gaifman graph.
For many of the following results it is therefore
sufficient to consider graphs.
\begin{lem}
Let\/ $\frakA$~be a structure.
A family $(U_v)_{v \in T}$ is a tree decomposition of\/~$\frakA$
if and only if it is a tree decomposition of\/ $\GF(\frakA)$.
\end{lem}
\begin{Proof}
$(\Rightarrow)$ is immediate. $(\Leftarrow)$
follows from the fact that every tree decomposition of a clique
has one component~$U_v$ containing the whole clique.
This implies that, for every clique~$C$ in $\GF(\frakA)$ induced
by some tuple $\bar c \in R_i$, there is some
vertex $v \in T$ with $C \subseteq U_v$.
Hence, every tuple $\bar c \in R_i$ is contained in some component~$U_v$.
\end{Proof}

In order to separate the higher classes of the hierarchy,
we shall employ two deep results of Robertson and Seymour about excluded minors.
\begin{thm}[Excluded Tree Theorem \cite{RobertsonSeymour83}]%
\label{Thm: excluded tree theorem}
For each tree~$\frakT$, there exists a number $k \in \bbN$ such that
\begin{align*}
  \frakT \notin \Min(\frakG) \qtextq{implies} \pwd(\frakG) < k\,,
  \qquad\text{for every graph } \frakG\,.
\end{align*}
\end{thm}
\begin{thm}[Excluded Grid Theorem \cite{RobertsonSeymour86}]%
\label{Thm: excluded grid theorem}
For each planar graph~$\frakE$, there exists a number $k \in \bbN$ such that
\begin{align*}
  \frakE \notin \Min(\frakG) \qtextq{implies} \twd(\frakG) < k\,,
  \qquad\text{for every graph } \frakG\,.
\end{align*}
\end{thm}
\begin{cor}
\itm{(a)} A class of graphs has bounded path-width if and only if
it excludes some tree as a minor.

\itm{(b)} A class of graphs has bounded tree-width if and only if
it excludes some planar graph as a minor.
\end{cor}
We also need a variant of these theorems for $n$-depth tree-width.
The next lemma contains the main technical argument.
\begin{lem}
Suppose that\/ $\frakG$~is a graph that does not contain a path of length~$l$.
Then $\frakG$~has a tree decomposition of height at most~$l$ and
width at most $l-1$.
\end{lem}
\begin{Proof}
Let $\langle T, {\preceq}\rangle$ be a depth-first spanning order-tree of~$\frakG$, i.e.,
a spanning tree such that, for every edge $(u,v)$ of~$\frakG$
we have $u \preceq v$ or $v \preceq u$
(for details see, e.g., \cite{Diestel06} where such spanning trees are called \emph{normal}).
We define a tree decomposition $(U_v)_{v \in T}$ of~$\frakG$ by setting
\begin{align*}
  U_v := \set{ u \in T }{ u \preceq v }\,.
\end{align*}
Since $T$~is depth-first, it follows that every edge $(u,v)$ of~$\frakG$
is contained in some component~$U_w$ where $w$~is the maximum of $u$~and~$v$.

The height of the tree~$T$ can be at most~$l$
since $\frakG$~contains no path of length~$l$. Furthermore, we have
$\abs{U_v} = \abs{v}+1 \leq l$.
Hence, the width of the tree decomposition is at most $l-1$.
\end{Proof}

\begin{thm}[Excluded Path Theorem]\label{Thm: excluded path theorem}
For each path~$\frakP$,
there exist numbers $n,k \in \bbN$ such that
\begin{align*}
  \frakP \notin \Min(\frakG) \qtextq{implies} \twd_n(\frakG) < k\,,
  \qquad\text{for every graph } \frakG\,.
\end{align*}
\end{thm}
\begin{Proof}
Suppose that $\frakP \notin \Min(\frakG)$ and let $l$~be the length of~$\frakP$.
Then the preceding lemma implies that $\twd_l(\frakG) < l$.
\end{Proof}
\begin{cor}
\itm{(a)} A class of graphs has bounded $n$-depth tree-width, for some~$n$, if and only if
it excludes some path as a minor \textup(equivalently, as a subgraph\textup).

\itm{(b)} A class of graphs has bounded tree-depth if and only if
it excludes some path as a minor \textup(equivalently, as a subgraph\textup).
\end{cor}

We can also compute a bound on the $n$-depth tree-width
in terms of the $(n+1)$-depth tree-width.
It will be needed in the proof of
Theorem~\ref{Thm: covering for lower part} below.

We say that the tree $\langle S,{\preceq}\rangle$ can be \emph{embedded}
into a tree~$\langle T,{\preceq}\rangle$
if there exists an order-preserving injective mapping
$\langle S, {\preceq}\rangle \to \langle T, {\preceq}\rangle$,
i.e., if $\langle S, {\preceq}\rangle$, regarded as relational structure,
is isomorphic to an induced substructure of $\langle T, {\preceq}\rangle$.
For instance, we have an embedding
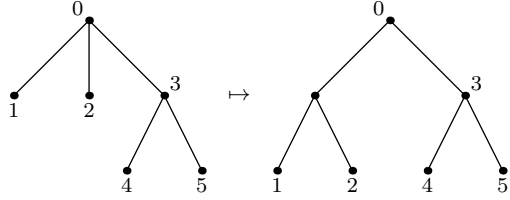
\begin{center}
%
%
%
%
%
%
%
\includegraphics{Transductions-final.1.ps}
\end{center}
as indicated by the labels.
If $S$~can be embedded in~$T$ then
$S$~is isomorphic to a minor of~$T$, when we consider
$S$~and~$T$ as graphs.

\begin{defi}
Let $D = (U_v)_{v \in T}$ be a tree decomposition and let $F$~be a set of edges of
(the successor-tree corresponding to)~$T$.
The tree decomposition $D/F$ obtained by \emph{contracting} the edges in~$F$ is
\begin{align*}
  D/F := (U'_{[v]})_{[v] \in T/F}\,,
\end{align*}
where
$U'_{[v]} := \bigcup_{u \in [v]} U_u$\,.
\end{defi}

\begin{lem}\label{Lem: bound on twd-n}
Let $\frakG$~be a graph and let $D := (U_v)_{v \in T}$ be a
tree decomposition of~$\frakG$ of width~$k$ and height at most $n+1$.
If $m \in \bbN$ is some number such that
the tree $m^{<n+1}$ cannot be embedded into~$T$,
then $\twd_n(\frakG) < m(k+1)$.
\end{lem}
\begin{Proof}
We construct a tree decomposition~$D'$ of
height at most~$n$ and width at most $m(k+1)-1$ as follows.
Let $P \subseteq T$ be
the minimal (w.r.t.~$\subseteq$) set of vertices that contains
\begin{itemize}
\item every leaf of~$T$ at level~$n$ and
\item every vertex that has at least $m$~successors in~$P$.
\end{itemize}
Since $m^{<n+1}$ cannot be embedded into~$T$ it follows that
$P$~does not contain the root of~$T$.
Let $F$~be the set of all edges of~$T$ linking
a vertex in $T \smallsetminus P$ to a vertex in~$P$. By definition of~$P$
it follows that (i)~every vertex of~$T$ has less than~$m$ $F$-successors\?;
(ii)~every path of~$T$ from the root to some leaf on level~$n$
contains at least one edge from~$F$\?; and (iii)~no such path contains two consecutive edges
from~$F$.

The decomposition $D' := D/F$
obtained by contracting all edges in~$F$
has width at most
\begin{align*}
  k+1 + (m-1)(k+1) - 1 < m(k+1)\,.
\end{align*}
Furthermore, the height of the underlying tree is at most~$n$.
\end{Proof}
\begin{cor}\label{Cor: bound on twd-n}
Let $\calC$~be a class such that $k := \twd_{n+1}(\calC) < \infty$
and let $m \in \bbN$.
If every structure $\frakA \in \calC$
has a tree decomposition $(U_v)_{v \in T}$ of width~$k$ and height at most $n+1$
such that the tree $m^{<n+1}$ cannot be embedded into~$T$,
then $\twd_n(\calC) < m(k+1) < \infty$.
\end{cor}

We conclude this section with a lemma that will be useful
when constructing transductions~$\tau_n$ that transform a structure
into their tree decompositions of height~$n$.
Our construction works for all tree decompositions that are \emph{strict}
in the following sense.
\begin{defi}\label{Def: strict tree decomposition}
Let $(U_v)_{v \in T}$ be a tree decomposition of a structure~$\frakA$.
\begin{enuma}
\item We define a function $\mu : A \to T$ by
  \begin{align*}
    \mu(a) := \min {\set{ v \in T }{ a \in U_v }}\,.
  \end{align*}
  Note that $\mu(a)$~is well-defined since, by the definition of a tree decomposition,
  there is at least one $v \in T$ such that $a \in U_v$.
\item For $v \in T$, we set
  \begin{align*}
    U_{\Aboveseg v} := \bigcup_{u \succeq v} U_u \smallsetminus \bigcup_{u \prec v} U_u\,.
  \end{align*}
\item The tree decomposition $(U_v)_v$ is \emph{strict} if, for every $v \in T$,
  \begin{itemize}
  \item $U_v \cap \mu(A) \neq \emptyset$ (equivalently, $U_v \smallsetminus U_u \neq \emptyset$,
    where $u$~is the predecessor of~$v$) and
  \item if $v$~is not the root of~$T$,
    then the subgraph of $\GF(\frakA)$ induced by the set $U_{\Aboveseg v}$ is connected.
  \end{itemize}
\end{enuma}
\end{defi}

We conclude this section by a result implying that, for our purposes,
it will be sufficient to consider only strict tree decompositions.
\begin{lem}\label{Lem: existence of strict tree decompositions}
Let $\frakG$~be a graph.
For every tree decomposition $(U_v)_{v \in T}$ of~$\frakG$,
there exists a strict tree decomposition $(U'_v)_{v \in T'}$ of~$\frakG$
whose width and height are at most those of $(U_v)_{v \in T}$.
\end{lem}
\begin{Proof}
By induction on $n \in \bbN$, we will construct a sequence $(U^n_v)_{v \in T_n}$
of tree decompositions such that $U^n_{\Aboveseg v}$~is connected,
for every $v \in T_n$ with level $0 < \abs{v} \leq n$.
(Recall that $\abs{v}$ denotes the level of~$v$, and the root is the only vertex of level~$0$.)
Furthermore, the restriction of~$T_n$ to the set of vertices of level at most~$n$
will coincide with the corresponding restriction of~$T_{n+1}$,
and we have $U^{n+1}_v = U^n_v$, for all $v \in T_{n+1}$ with level $\abs{v} \leq n$.
It will follow that the sequence has a limit $(U^\omega_v)_{v \in T_\omega}$
where
\begin{align*}
  T_\omega := \bigcup_{n \in \bbN} {\set{ v \in T_n }{ \abs{v} \leq n }}
  \qtextq{and}
  U^\omega_v := U^{\abs{v}}_v\,.
\end{align*}

We start the construction with $T_0 := T$ and $U^0_v := U_v$.
Suppose that we have already defined $(U^n_v)_{v \in T_n}$.
For every vertex $v \in T_n$ of level $\abs{v} = n+1$
we modify the tree decomposition as follows.
Let $C_0,\dots,C_{m-1}$ be an enumeration of the connected components
of~$U^n_{\Aboveseg v}$. We replace in~$T_n$ the subtree rooted at~$v$
by $m$~copies $S_0,\dots,S_{m-1}$ of the subtree, all attached
to the predecessor of~$v$. For $u \in S_i$ we define
$U^{n+1}_u := U^n_u \cap C_i$.
We can do these modifications for all vertices of level $n+1$ simultaneously.
Let $(U^{n+1}_v)_{v \in T_{n+1}}$ be the resulting tree decomposition.

The limit $(U^\omega_v)_{v \in T_\omega}$ of this sequence satisfies
the connectedness requirement of a strict tree decomposition.
To also satisfy the other condition we proceed as follows.
Let $F$~be the set of all edges $(u,v)$ of~$T_\omega$ such that
$U_v \cap \mu(V) = \emptyset$.
(Note that this implies $U_v \subseteq U_u$.)
We construct the tree decomposition $(U'_v)_{v \in T'}$
by contracting all edges in~$F$.
The details and the remaining verifications are left to the reader.
\end{Proof}

\section{Tree decompositions and transductions}   
\label{Sect: tree decompositions and transductions}

In this section we relate the material presented in the preceding one
to monadic second-order transductions.
Let us start by showing that there is a transduction computing the minors
of a graph.
\begin{lem}[\cite{Courcelle95}]\label{Lem: minors interpretable}
There exists a transduction~$\tau$ such that
$\tau(\frakG_\IN) = \Min(\frakG)$, for every graph~$\frakG$.
\end{lem}
\begin{Proof}
A minor~$\frakH$ of~$\frakG$ is obtained by
deleting vertices, deleting edges, and contracting edges.
Hence, we can encode~$\frakH$ by four sets\?:
the set of vertices we delete, the set of edges
we delete, the set of edges we contract, and a set of vertices
containing one representative of each contracted subgraph of~$\frakG$
(these vertices serve as vertices of the resulting graph~$\frakH$).
With the help of these parameters we can define~$\frakH$ inside of~$\frakG_\IN$
by $\MSO$-formulae.
\end{Proof}

There is a close relationship between tree decompositions and
transductions.
\begin{lem}\label{Lem: tree decomposable implies interpretable}
For every signature~$\Sigma$ and every number $k \in \bbN$,
there exists a transduction $\tau_k : \TREE_0 \to \STR_\IN[\Sigma]$
that maps an order-tree~$T$ to the class of all incidence structures~$\frakA_\IN$
such that the corresponding $\Sigma$-structure~$\frakA$ has
a tree decomposition of width at most~$k$ with underlying tree~$T$.
\end{lem}
\begin{Proof}
Suppose that $\frakA$~is a structure which has a tree decomposition
$(U_v)_{v \in T}$ of width~$k$.
We prove that $\frakA$~can be defined from a colouring of~$T$
where the number of colours depends only on $\Sigma$~and~$k$.

Let $\frakC_0,\dots,\frakC_{m-1}$ be an enumeration
of all $\Sigma$-structures whose domain is a subset of $[k+1]$.
For each $v \in T$, let $\frakU_v$~be the substructure of~$\frakA$
induced by~$U_v$.
It follows that, for every $v \in T$, we can find some index $\lambda(v)$
such that $\frakU_v \cong \frakC_{\lambda(v)}$.
Let $\pi_v : \frakU_v \to \frakC_{\lambda(v)}$ be the corresponding isomorphism.

Furthermore, we associate with each edge $(u,v)$ of~$T$ the binary relation
\begin{align*}
  R(u,v) := \set{ (\pi_u(a), \pi_v(a)) }{ a \in U_u \cap U_v } \subseteq [k+1] \times [k+1]\,.
\end{align*}

We can recover~$\frakA$ from~$T$ with the help of the vertex colouring~$\lambda$
and the edge colouring~$R$.
We form the disjoint union of all structures $(\frakC_{\lambda(v)})_\IN$, for $v \in T$,
and we identify two elements $i \in C_{\lambda(u)}$ and $j \in C_{\lambda(v)}$
if $(u,v)$ is an edge of~$T$ such that $(i,j) \in R(u,v)$.
This can be performed by an $n$-copying $\MSO$-transduction
where $n$~is the maximal size of the structures $(\frakC_i)_\IN$, $i < m$.
\end{Proof}

We have just seen that we can map a class of trees to a class
of structures with these trees as tree decompositions.
Conversely, if we only consider strict tree decompositions,
we can define a transduction mapping a class of structures
to the corresponding class of trees.
Recall the function $\mu : A \to T$ from Definition~\ref{Def: strict tree decomposition}.
that assigns to an element $a \in A$ the minimal index $v \in T$ such that
$a \in U_v$.
\begin{prop}\label{Prop: strict tree decompositions definable}
For each number $n \in \bbN$, there exists an $\MSO$-formula
$\varphi_n(x,y;\bar Z)$ such that, for every strict tree
decomposition $D = (U_v)_{v \in T}$ of a graph~$\frakG$
of height at most~$n$, there are sets
$L_0,\dots,L_{n-1} \subseteq V$ such that
\begin{align*}
  \frakG \models \varphi_n(a,b;\bar L) \quad\iff\quad \mu(a) \leq \mu(b)\,.
\end{align*}
\end{prop}
\begin{Proof}
Given~$D$ we use the sets
\begin{align*}
  L_i := \set{ a \in V }{ \abs{\mu(a)} = i }
\end{align*}
of all elements that first appear at level~$i$ of the tree.
In particular, $L_0 = U_{\emptyseq}$ is the root component of the tree decomposition.
For $k < n$, let $\frakG_{\geq k}$~be the subgraph of~$\frakG$
induced by $L_k \cup\dots\cup L_{n-1}$.
For $a \in L_i$ and $b \in L_j$ we define
\begin{align*}
  a \preceq b
  \quad\iff\quad
  &i \leq j \text{ and $a$, $b$ belong to the same connected component of } \frakG_{\geq i}\,,\\
\prefixtext{and}
  a \sim b \quad\iff\quad &a \preceq b \text{ and } b \preceq a\,, \text{ or if } a,b \in L_0\,.
\end{align*}
Clearly, the relation~$\preceq$ is $\MSO$-definable with the help of the
parameters~$\bar L$. We claim that, for $a,b \notin L_0$, we have
\begin{align*}
  a \preceq b \quad\iff\quad \mu(a) \leq \mu(b)\,.
\end{align*}

$(\Leftarrow)$ Suppose that $\mu(a) \leq \mu(b)$.
Then $b \in U_{\Aboveseg\mu(a)}$. Furthermore,
$U_{\Aboveseg\mu(a)}$ is connected since $D$~is strict.
Hence, $U_{\Aboveseg\mu(a)}$ is a connected component of $\frakG_{\geq i}$
containing both $a$~and~$b$. Since $\abs{\mu(a)} \leq \abs{\mu(b)}$
it follows that $a \preceq b$.

$(\Rightarrow)$ Suppose that $a \preceq b$. Then there exists
an undirected path~$\pi$ in~$\frakG_{\geq\abs{\mu(a)}}$ connecting $a$~and~$b$.
Since $U_{\Aboveseg u} \cap U_{\Aboveseg v} = \emptyset$, for all $u \neq v$
such that $\abs{u} = \abs{v}$, it follows that $\pi$~is contained in some $U_{\Aboveseg v}$
such that $\abs{v} = \abs{\mu(a)}$. Since $a$~is a vertex of~$\pi$ we must have
$v = \mu(a)$. Furthermore, $b \in U_{\Aboveseg\mu(a)}$ since $b$~is also a vertex of~$\pi$.
This implies that $\mu(a) \leq \mu(b)$.
\end{Proof}
\begin{thm}\label{Thm: strict tree decompositions can be computed by a transduction}
For each constant $n \in \bbN$, there exists a transduction~$\tau_n$
mapping a graph~$\frakG$ to the class of all (underlying trees of)
strict tree decompositions of~$\frakG$ of height at most~$n$.
\end{thm}
\begin{Proof}
Let $D = (U_v)_{v \in T}$ be a strict tree decomposition of~$\frakG$, and
let $\varphi_n(x,y;\bar L)$ be the formula
from Proposition~\ref{Prop: strict tree decompositions definable}
with parameters $L_0,\dots,L_{n-1} \subseteq V$.
We can define the tree~$T$ underlying~$D$ as follows\?:
\begin{itemize}
\item Its root is any element of~$L_0 = U_{\emptyseq}$.
\item For the other vertices of~$T$, we choose one vertex in each $\sim$-class
  different from~$L_0$. Note that $\sim$~is definable with the help of~$\varphi_n$.
\item The ordering of~$T$ is defined by~$\varphi_n$.
\end{itemize}
Hence, we obtain a transduction with parameters~$\bar L$
that transforms a graph into a `candidate' tree decomposition.
Via a backwards translation we can write down a formula
stating that the candidate given by the parameters~$\bar L$
corresponds to an actual strict tree decomposition.
We omit the details which are standard for this type of construction.
\end{Proof}

In Lemma~\ref{Lem: tree decomposable implies interpretable}
we have seen how to obtain classes of bounded tree-width
from classes of trees. Conversely, it is the case that
every class obtained from a class of trees via a transduction
has a bounded tree-width.
\begin{thm}\label{Thm: interpretable implies tree decomposable}
For every transduction $\tau : \TREE_m \to \STR_\IN[\Sigma]$,
there exists a number $k \in \bbN$
such that, for each $m$-coloured tree~$\frakT$ with image $\frakA_\IN \in \tau(\frakT)$,
the structure $\frakA$~has a tree decomposition
of width at most~$k$ where the underlying tree is~$\frakT$.
\end{thm}
\begin{rem}
(a) Courcelle and Engelfriet~\cite{CourcelleEngelfriet95} have shown that
an incidence structure~$\frakA_\IN$ obtained via a transduction~$\tau$ from an
$m$-coloured tree~$\frakT$ has bounded tree-width.
Theorem~\ref{Thm: interpretable implies tree decomposable}
strengthens this result by proving that,
if $\frakA_\IN$~is the image of a tree~$\frakT$, then we can use the same
tree~$\frakT$ as the tree underlying a tree decomposition of the given width.

(b) Lapoire has announced in~\cite{Lapoire98} a result somewhat related to
Theorem~\ref{Thm: strict tree decompositions can be computed by a transduction}.
He claims that, for every $k \in \bbN$, there exists a transduction that transforms
a given graph~$\frakG$ of tree-width at most~$k$ to a coloured tree (like in the proof
of Lemma~\ref{Lem: tree decomposable implies interpretable})
that encodes some tree decomposition of~$\frakG$ of width at most~$k$.
Our result is less ambitious in the sense that we only consider tree decompositions
of a fixed height. This enables us to give a precise description of
which tree decompositions (the strict ones) our transduction returns.
Note that one can show that, for $k \geq 2$,
there is no such transduction that would return \emph{all} tree decompositions
of~$\frakG$ of width at most~$k$.
\end{rem}

We split the proof into several lemmas.
As a technical tool we introduce a second kind of hierarchical decompositions
of structures and a corresponding notion of width.
To simplify the definition we will only consider incidence structures.
\begin{defi}
Let $\frakA_\IN = \langle A \cup E, \bar P,\IN_0,\dots\rangle$ be an incidence
structure.

(a)
A \emph{partition refinement} of~$\frakA_\IN$ is a family $\Pi = (W_v, {\approx_v})_{v \in T}$
of pairs consisting of a subset $W_v \subseteq A \cup E$ and an equivalence relation~$\approx_v$
on~$W_v$ with the following properties\?:
\begin{itemize}
\item The index set~$T$ is a tree.
\item For the root~$\emptyseq$, we have $W_{\emptyseq} = A \cup E$
\item For every internal vertex (i.e., non-leaf) $u \in T$ with successors $v_0,\dots,v_{n-1}$,
  the sets $W_{v_0},\dots,W_{v_{n-1}}$ form a partition of~$W_u$.
\item $\abs{W_u} = 1$, for every leaf $u \in T$.
\item $x \approx_v y$ and $u \preceq v$ implies $x \approx_u y$.
\item If $u$~is an internal vertex of~$T$, $v,w$ successors of~$u$, not necessarily distinct,
  and $x \in W_v$, $y \in W_w$ elements, then
  $x \approx_u y$ implies either
  \begin{alignat*}{-1}
    &x,y \in A \text{ and, } &&\text{for every }  e \in E \smallsetminus (W_v \cup W_w) \text{ and every } i\,, \\
    &&& (x,e) \in \IN_i \Leftrightarrow (y,e) \in \IN_i
  \prefixtext{or}
    &x,y \in E \text{ and, } &&\text{for every } a \in A \smallsetminus (W_v \cup W_w) \text{ and every } i\,, \\
    &&& (a,x) \in \IN_i \Leftrightarrow (a,y) \in \IN_i\,.
  \end{alignat*}
\end{itemize}
Note that it follows that, for every element $x \in A \cup E$, there is some leaf $u \in T$
such that $W_u = \{x\}$.

(b)
The \emph{width} of a partition refinement $\Pi = (W_v, {\approx_v})_{v \in T}$ is the maximum
number of equivalence classes realised in some component~$W_v$\?:
\begin{align*}
  \WD(\Pi) := \max_{v \in T} {\abs{W_v/{\approx_v}}}\,.
\end{align*}
The \emph{partition-width} of the structure~$\frakA_\IN$ is the minimal width
of a partition refinement of~$\frakA_\IN$.
\end{defi}
The notion of a partition refinement and of partition-width
are adaptations of definitions from \cite{Blumensath06,Blumensath03b}.
Up to a factor of~$2$, the partition-width of an incidence structure and its clique-width
coincide.
\begin{exa}
Let $\frakA = (A,R)$ be a structure with domain $A = \{a,b,c,d,e\}$
and a ternary relation
\begin{align*}
  R = \{\underbrace{(a,b,c)}_x,\underbrace{(a,b,d)}_y,\underbrace{(a,b,e)}_z\}\,.
\end{align*}
Its incidence structure is $\frakA_\IN = \langle A \cup E, P_R, \IN_0,\IN_1,\IN_2\rangle$
with $E = \{x,y,z\}$.
We obtain a partition refinement
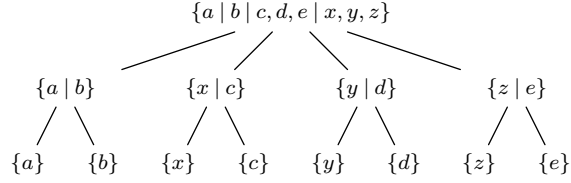
\begin{center}
%
%
%
%
\includegraphics{Transductions-final.2.ps}
\end{center}
where we have indicated the partition into $\approx_v$-classes by vertical bars.
This partition refinement has width~$4$.
\end{exa}

\begin{lem}\label{Lem: tree decomposition bounded by partition width}
For every partition refinement $\Pi = (W_v, {\approx_v})_{v \in T}$
of an incidence structure\/
$\frakA_\IN = \langle A \cup E,\bar P,\IN_0,\dots,\IN_{r-1}\rangle$,
there exists a tree decomposition $D = (U_v)_{v \in T}$
of\/~$\frakA$ with the same underlying tree~$T$ such that
\begin{align*}
  \WD(D) < (r+3)\cdot\WD(\Pi)\,.
\end{align*}
\end{lem}
\begin{Proof}
Let $l : A \cup E \to \Lf(T)$ be the function assigning to every $x \in A \cup E$
the unique leaf $l(x)$ of~$T$ such that $W_{l(x)} = \{x\}$.
We claim that the desired tree decomposition $(U_u)_{u \in T}$ of~$\frakA$
is given by
\begin{align*}
  U_u &:= B_u \cup C_u \cup D_u
\intertext{where}
  B_u &:= \set{ v \in A }
              { u \preceq l(v) \text{ and } (v,e) \in \IN_i \text{ for some } i < r \text{ and } e \in E \text{ with }
                u \npreceq l(e) }\,, \\
  C_u &:= \set{ v \in A }
              { u \npreceq l(v) \text{ and } (v,e) \in \IN_i \text{ for some } i < r \text{ and } e \in E \text{ with }
                u \preceq l(e) }\,, \\
  D_u &:= \set{ v \in A }
              { (v,e) \in \IN_i \text{ for some } i < r \text{ and } e \in E \text{ with } 
                l(v) \sqcap l(e) = u }\,.
\end{align*}
Note that the connectedness condition holds since $(v,e) \in \IN_i$ implies that
$v$~belongs to precisely those components~$U_u$ such that~$u$ lies on the path
from~$l(v)$ to~$l(e)$.

It remains to prove that $\abs{U_u} \leq (r+3)\cdot\WD \Pi$.
If $u = l(\bar c)$, for some $\bar c \in E$, then $U_u = C_u$
consists of the components of~$\bar c$.
Hence, $\abs{U_u} = \abs{\bar c} \leq r$.
Therefore, we may assume that $u \notin l[E]$. Let
\begin{align*}
  [x]_u := \set{ y \in W_u }{ y \approx_u x }
\end{align*}
denote the $\approx_u$-class of~$x$.
We prove the following bounds.
\begin{enumr}
\item $\abs{[x]_u} = 1$, for all $x \in B_u$.
\item $\abs{[x]_u \cap U_u} \leq 2$, for all $x \in D_u$.
\item $\abs{C_u} \leq r \cdot \abs{W_u/{\approx_u}}$.
\end{enumr}
Since $B_u,D_u \subseteq W_u$ it then follows
that $\abs{U_u} = \abs{B_u \cup C_u \cup D_u} \leq (r+3)\cdot\abs{W_u/{\approx_u}}$.

(1) Let $x \in B_u$. There is some tuple $e \in E$ and some index~$i$ such that $(x,e) \in \IN_i$
and $u \npreceq l(e)$.
We have $(y,e) \in \IN_i$, for every $y \in W_u$ such that $y \approx_u x$.
Since $x$~is the only such element it follows that $[x]_u = \{x\}$.

(2) Let $x \in D_u$. There is some tuple $e \in E$ and some~$i$ such that
$(x,e) \in \IN_i$ and $u = l(x) \sqcap l(e)$.
Let $v$~be the successor of~$u$ such that $v \preceq l(e)$.
We have $(y,e) \in \IN_i$, for all $y \in W_u \smallsetminus W_v$ such that $y \approx_u x$.
Hence, $[x]_u \smallsetminus W_v = \{x\}$.

Suppose that there is some element $y \in [x]_u \cap W_v \cap U_u$.
By definition of~$U_u$ there is some tuple~$f \in E$ and some~$j$ such that $(y,f) \in \IN_j$
and $l(y) \sqcap l(f) \preceq u$.
As above it follows that $[x]_u \cap W_v = \{y\}$.
Consequently, we have $\abs{[x]_u \cap U_u} \leq 2$.

(3) Let $x \in C_u$ and consider some tuple $e \in E$ such that $(x,e) \in \IN_i$ and $u \preceq l(e)$.
Set
\begin{align*}
  I_u(e) := \set{ z \in A }{ (z,e) \in \IN_i \text{ for some } i \text{ and } u \npreceq l(z) }\,.
\end{align*}
For $e,f \in E \cap W_u$, it follows that
\begin{align*}
  e \approx_u f
  \qtextq{implies}
  I_u(e) = I_u(f)\,.
\end{align*}
Furthermore, we obviously
have $\abs{I_u(e)} \leq \abs{e} \leq r$.
It follows that $C_u$~contains at most $r\cdot \abs{W_u/{\approx_u}}$ vertices.
\end{Proof}

\begin{lem}\label{Lem: partition width of simple tree interpretations}
Let $\tau : \TREE_m \to \STR_\IN[\Sigma]$ be a basic $\MSO$-transduction such that,
for every $m$-coloured order-tree~$\frakT$ with image $\frakA_\IN \in \tau(\frakT)$, we have
\begin{align*}
  A \cup E = \Lf(\frakT)
  \qtextq{and}
  A \cap E = \emptyset\,.
\end{align*}
Then there exists a number $n \in \bbN$ such that,
for every order-tree~$\frakT$, we can find a partition refinement $(W_v, {\approx_v})_{v \in T}$
of $\tau(\frakT)$ of width at most~$n$.
\end{lem}
\begin{Proof}
Let $\langle\chi,\delta(x),(\varphi_{P_R}(x))_R,(\varphi_{\IN_i}(x,y))_{i < r}\rangle$
be the definition scheme of~$\tau$, and
let $h$~be the maximal rank of these formulae.

Given~$\frakT$ we define the desired partition refinement $\Pi = (W_u, {\approx_u})_{u \in T}$
by setting
\begin{align*}
  &W_u := \set{ x \in \Lf(T) }{ u \preceq x }\,, \\[1ex]
\prefixtext{and}
  &x \approx_u y \quad\defiff\quad
  \begin{aligned}[t]
    &\MTh_h(\frakT_v,x) = \MTh_h(\frakT_w,y)\,, \\
    &\text{where $v,w$ are the successors of $u$ with }
     x \in W_v \text{ and } y \in W_w\,.
  \end{aligned}
\end{align*}
(If $u$~is a leaf of~$T$ then $W_u = \{x\}$ and we take the equality relation for~$\approx_u$.)
Note that the index of~$\approx_v$ is finite and
that it only depends on~$h$ and not on the input tree~$\frakT$.

It remains to show that $\Pi$~is actually a partition refinement.
First, let us prove that $x \approx_v y$ and $u \preceq v$ implies $x \approx_u y$.
It is sufficient to consider the case that $u$~is the predecessor of~$v$.
Then the general case follows by induction.
Hence, suppose that $v$~is a successor of~$u$, that $w,w'$ are
successors of~$v$, and that $x,y$ are leaves with $w \preceq x$ and $w' \preceq y$
such that $x \approx_v y$. Then we have
\begin{align*}
  \MTh_h(\frakT_w,x) = \MTh_h(\frakT_{w'},y)\,,
\end{align*}
which, by Lemma~\ref{Lem: composition for trees}, implies that
\begin{align*}
  \MTh_h(\frakT_v,x) = \MTh_h(\frakT_v,y)\,.
\end{align*}
Consequently, we have $x \approx_u y$.

We also have to show that the incidence relation is invariant under~$\approx_u$.
Let $v,w$~be successors of~$u$ and suppose that $x,y$ are leaves with
$v \preceq x$ and $w \preceq y$ such that $x \approx_u y$.
We distinguish two cases.

Suppose that $x,y \in A$ and let $e \in E \smallsetminus (W_v \cup W_w)$ be an edge.
Since 
\begin{align*}
  \MTh_h(\frakT_v,x) = \MTh_h(\frakT_w,y)\,,
\end{align*}
it follows that
\begin{align*}
  \frakT \models \varphi_{\IN_i}(x,e) \quad\iff\quad
  \frakT \models \varphi_{\IN_i}(y,e)\,.
\end{align*}
Hence, $(x,e) \in \IN_i$ iff $(y,e) \in \IN_i$.

Now, suppose that $x,y \in E$ and let $z \in A \smallsetminus (W_v \cup W_w)$ be an element.
Since 
\begin{align*}
  \MTh_h(\frakT_v,x) = \MTh_h(\frakT_w,y)\,,
\end{align*}
it follows that
\begin{align*}
  \frakT \models \varphi_{\IN_i}(z,x) \quad\iff\quad
  \frakT \models \varphi_{\IN_i}(z,y)\,.
\end{align*}
Hence, $(z,x) \in \IN_i$ iff $(z,y) \in \IN_i$.
\end{Proof}

\begin{Proof}[Proof of Theorem~\ref{Thm: interpretable implies tree decomposable}]
(1) First, suppose that $\tau : \TREE_m \to \STR_\IN[\Sigma]$
is a basic $\MSO$-transduction such that,
for every $m$-coloured order-tree $\frakT \in \dom(\tau)$ with image $\frakA_\IN \in \tau(\frakT)$,
we have
\begin{align*}
  A \cup E = \Lf(\frakT)
  \qtextq{and}
  A \cap E = \emptyset\,.
\end{align*}
It follows by Lemma~\ref{Lem: partition width of simple tree interpretations} that
there is a number $w \in \bbN$ such that,
for every tree~$\frakT$, we can find a partition refinement of
$\frakA_\IN \in \tau(\frakT)$
with underlying tree~$\frakT$ whose width is at most~$w$.
By Lemma~\ref{Lem: tree decomposition bounded by partition width}
it follows that $\frakA$~has a tree decomposition $(U_v)_{v \in T}$
with underlying tree~$\frakT$ and whose width is less than $k := w(r+3)$.

(2) If $\tau$~is a basic $\MSO$-transduction such that
\begin{align*}
  A \cup E \subseteq \Lf(\frakT)
  \qtextq{and}
  A \cap E = \emptyset\,,
  \quad\text{for all } \frakA_\IN \in \tau(\frakT) \text{ with } \frakT \in \dom(\tau)\,,
\end{align*}
then we can argue similarly.
Let $\tau'$~be the $\MSO$-transduction mapping~$\frakT$ to the structure
obtained from $\tau(\frakT)$ by adding one isolated element for every leaf of~$\frakT$
that does not correspond to an element of~$\tau(\frakT)$.
Then $\tau'$~is of the form considered in~(1) and we obtain a tree decomposition~$(U_v)_{v \in T}$
of $\tau'(\frakT)$. Deleting from every component~$U_v$
all elements not in~$\tau(\frakT)$ we obtain the desired tree decomposition of~$\tau(\frakT)$.

(3) Suppose that $\tau$~is a non-copying $\MSO$-transduction as in~(2) but with~$p$ parameters.
We can regard~$\tau$ as a basic $\MSO$-transduction $\TREE_{m+p} \to \STR_\IN[\Sigma]$.
By~(2) it follows that,
for every value of the parameters~$\bar P$,
the structure $\tau(\frakT,\bar P)$ has a tree decomposition of the required form.

(4) Finally, consider the general case. Suppose that $\tau$~is $l$-copying.
Given~$\frakT$ let $\frakT^+$~be the tree obtained from~$\frakT$
by adding~$l$ new successors to every vertex of~$\frakT$.
Formally, suppose that $T \subseteq D^{<\omega}$, for some finite set~$D$.
W.l.o.g.\ we may assume that $D \cap [l] = \emptyset$.
We define the domain $T^+ \subseteq (D \cup [l])^{<\omega}$ of~$\frakT^+$ by
\begin{align*}
  T^+ := T \cup (T \times [l])\,.
\end{align*}
Furthermore, we add new colour predicates
\begin{align*}
  S_i := T \times \{i\}\,, \quad\text{for } i \in [l]\,.
\end{align*}
Note that every element of $\tau(\frakT)$ is of the form $\langle v,i\rangle$
where $i \in [l]$ and $v \in T$. Hence, each such element corresponds to a leaf
$vi \in T \times [l] \subseteq T^+$.
Using the parameters~$\bar S$ we can construct a
basic $\MSO$-transduction $\tau^+ : \TREE_{m+l} \to \STR_\IN[\Sigma]$ 
satisfying the conditions of~(3) such that
$\tau^+(\frakT^+) = \tau(\frakT)$.
By~(3), we obtain a tree decomposition $H^+ = (U^+_v)_{v \in T^+}$ of $\tau^+(\frakT^+) = \tau(\frakT)$.
Let $H = (U_v)_{v \in T}$ be the tree decomposition obtained from~$H^+$ by
contracting every edge leading to a leaf in $T^+ \smallsetminus T$.
Then we have
\begin{align*}
  \WD(H) + 1 &\leq (l+1)(\WD(H^+)+1)
\prefixtext{and}
  \WD(H^+) &\leq w\,, \quad\text{for some } w \in \bbN \text{ independent of } \frakT\,.
\end{align*}
\upqed
\end{Proof}

\section{The transduction hierarchy}   
\label{Sect: hierarchy}

The focus of our investigation lies on the following preorder on classes of structures
which compares their `encoding powers' with respect to $\MSO$-trans\-duc\-tions.
Our main result is a complete description of the hierarchy induced by this preorder.
It will be given in Theorem~\ref{Thm: hierarchy}.
\begin{defi}
Let $\calC,\calK \subseteq \STR$. We define the following relations.
\begin{enuma}
\item $\calC \sqsubseteq \calK$ if there exists a transduction~$\tau$
  such that $\calC \subseteq \tau(\calK)$.
\item $\calC \sqsubset \calK$ if $\calC \sqsubseteq \calK$ and $\calK \nsqsubseteq \calC$.
\item $\calC \equiv \calK$ if $\calC \sqsubseteq \calK$ and $\calK \sqsubseteq \calC$.
\item $\calC \Sqsubset \calK$ if $\calC \sqsubset \calK$ and there is no class~$\calD$
  with $\calC \sqsubset \calD \sqsubset \calK$.
\item $\calC \sqsubseteq_\IN \calK$ if $\calC_\IN \sqsubseteq \calK_\IN$.
\item The relations $\sqsubset_\IN$, $\equiv_\IN$, and $\Sqsubset_\IN$ are defined analogously
  to $\sqsubset$, $\equiv$, $\Sqsubset$ by replacing~$\sqsubseteq$ everywhere by~$\sqsubseteq_\IN$.
\end{enuma}
The \emph{transduction hierarchy} is the hierarchy of classes $\calC \subseteq \STR$
induced by the relation~$\sqsubseteq_\IN$.
\end{defi}

As the class of transductions is closed under composition, it follows that
the relation~$\sqsubseteq_\IN$ is a preorder, i.e., it is reflexive and transitive.
\begin{lem}
$\sqsubseteq_\IN$~is a preorder on $\PSet(\STR)$.
\end{lem}

\begin{defi}
We consider the following subclasses of $\STR[\{\mathrm{edg}\}]$.
(All trees below are considered to be successor-trees.)
\begin{enuma}
\item $\calT_n := \set{ m^{<n} }{ m \in \bbN}$ is the set of all complete $m$-ary trees of height~$n$.
\item $\calT_{\mathrm{bin}}$ is the class of all binary trees.
\item $\calT_\omega$ is the class of all trees.
\item $\calP$~is the class of all paths.
\item $\calG$~is the class of all rectangular grids.
\end{enuma}
\end{defi}

The following description of the transduction hierarchy is the main result of the present paper.
\begin{thm}\label{Thm: hierarchy}
We have the following hierarchy\?:
\begin{align*}
  \emptyset \Sqsubset_\IN \calT_0 \Sqsubset_\IN \calT_1 \Sqsubset_\IN\dots\Sqsubset_\IN \calT_n \Sqsubset_\IN\dots
  \sqsubset_\IN \calP \Sqsubset_\IN \calT_\omega \equiv_\IN \calT_{\mathrm{bin}} \Sqsubset_\IN \calG
\end{align*}
For every signature~$\Sigma$,
every class $\calC \subseteq \STR[\Sigma]$ is $\equiv_\IN$-equivalent to some class in this hierarchy.
\end{thm}
\begin{rem}
There is a lot of flexibility in the choice
of representatives for the various levels.
For instance, we could replace~$\calT_n$ by the class of \emph{all} trees of height at most~$n$,
$\calT_\omega$ by $\calT_{\mathrm{bin}}$ or $\set{ 2^{<n} }{ n \in \bbN }$,
and $\calG$~by the class of square grids.
\end{rem}

It is straightforward to show that the above classes form an increasing chain.
The hard part is to prove that the chain is strictly increasing and that
there are no further classes.
\begin{lem}
We have
\begin{align*}
  \emptyset \sqsubseteq_\IN \calT_0 \sqsubseteq_\IN \calT_1 \sqsubseteq_\IN\dots\sqsubseteq_\IN \calT_n \sqsubseteq_\IN\dots
  \sqsubseteq_\IN \calP \sqsubseteq_\IN \calT_\omega \sqsubseteq_\IN \calG\,.
\end{align*}
\end{lem}
\begin{Proof}
In the example before Lemma~\ref{Lem: transductions are comorphisms},
we have constructed transductions~$\tau_n$ such that
$\calT_n \subseteq \tau_n(\calP)$. Hence, $\calT_n \sqsubseteq_\IN \calP$.
The remaining assertions follow from the observation that,
by Lemma~\ref{Lem: minors interpretable},
$\calC \subseteq \Min(\calK)$ implies $\calC \sqsubseteq_\IN \calK$.
\end{Proof}

Let us collect some easy properties of the hierarchy.
Our first result states that $\calG$~is a representative of the top level
of the transduction hierarchy.
\begin{lem}\label{Lem: G is top}
$\STR[\Sigma] \sqsubseteq_\IN \calG$
\end{lem}
\begin{Proof}
Recall that the $m \times n$ grid is the undirected graph
$\frakG = \langle V,\edg\rangle$ with vertices
$V = [m] \times [n]$ and edge relation
\begin{align*}
  \edg = \bigset{ (\langle i,k\rangle,\langle j,l\rangle) }{ \abs{i-j} + \abs{k-l} = 1 }\,.
\end{align*}
Before encoding arbitrary structures in such grids we describe a transduction
mapping~$\frakG$ to its directed variant $\langle V, E_0, E_1\rangle$ where
\begin{align*}
  E_0 &:= \set{ (\langle i,k\rangle,\langle i+1,k\rangle) }{ i < m-1,\ k < n }\,, \\
\prefixtext{and}
  E_1 &:= \set{ (\langle i,k\rangle,\langle i,k+1\rangle) }{ i < m,\ k < n-1 }\,.
\end{align*}
This can be done with the help of the parameters $P_0,P_1,P_2,Q_0,Q_1,Q_2 \subseteq V$
where
\begin{align*}
  P_m &:= \set{ \langle i,k\rangle }{ i \equiv m \pmod 3 }\,, \\
\prefixtext{and}
  Q_m &:= \set{ \langle i,k\rangle }{ k \equiv m \pmod 3 }\,.
\end{align*}
Then
\begin{align*}
  E_0 &= \set{ (u,v) \in \edg }{ u \in P_i \text{ and } v \in P_j \text{ for some } i \equiv j-1 \pmod 3 }\,, \\
\prefixtext{and}
  E_1 &= \set{ (u,v) \in \edg }{ u \in Q_i \text{ and } v \in Q_j \text{ for some } i \equiv j-1 \pmod 3 }\,.
\end{align*}
It is easy to write down a formula checking that the parameters $P_m$~and~$Q_m$
are correctly chosen (see, e.g., \cite{Courcelle97}).

To show that $\STR[\Sigma] \sqsubseteq_\IN \calG$,
suppose that $\frakA \in \STR[\Sigma]$ is a structure with
$\frakA_\IN = \langle A \cup E, (P_R)_R, \IN_0,\dots,\IN_{r-1}\rangle$.
Fix enumerations $a_0,\dots,a_{m-1}$ of~$A$ and $e_0,\dots,e_{n-1}$ of~$E$.
By the above remarks, it is sufficient to encode~$\frakA_\IN$ in the directed
$m \times n$ grid.
Consider the following subsets of $[m] \times [n]$\?:
\begin{alignat*}{-1}
  A'   &:= [m] \times \{0\}\,,\qquad
  &P'_R &:= \set{ \langle 0,k\rangle }{ e_k \in P_R }\,, \\
  E'   &:= \{0\} \times [n]\,,\qquad
  &I'_l &:= \set{ \langle i,k\rangle }{ (a_i,e_k) \in \IN_l }\,.
\end{alignat*}
Then $\frakA_\IN$ can be recovered from~$\frakG$ by an $\MSO$-transduction
using these sets as parameters.
\end{Proof}

\begin{lem}
$\calT_\omega \equiv_\IN \calT_{\mathrm{bin}}$.
\end{lem}
\begin{Proof}
For one direction, note that $\calT_{\mathrm{bin}} \subseteq \calT_\omega$ implies
$\calT_{\mathrm{bin}} \sqsubseteq_\IN \calT_\omega$.
Conversely, each finite tree can be obtained as minor of a binary tree.
Therefore, we have $\calT_\omega \subseteq \Min(\calT_{\mathrm{bin}}) \sqsubseteq_\IN \calT_{\mathrm{bin}}$.
\end{Proof}

\begin{lem}
We have $\calC \equiv_\IN \calT_1$ if and only if $\calC$~is finite and contains at least one
nonempty structure.
\end{lem}

As indicated in the example before Lemma~\ref{Lem: transductions are comorphisms},
there exists a transduction mapping an incidence structure~$\frakA_\IN$ to the
incidence structure $\GF(\frakA)_\IN$ of the Gaifman graph of~$\frakA$.
\begin{lem}\label{Lem: Gaifman graph interpretable}
For every class~$\calC$ of structures, we have $\GF(\calC) \sqsubseteq_\IN \calC$\,.
\end{lem}
The next result is just a restatement of Lemma~\ref{Lem: minors interpretable}
in our current terminology.
\begin{lem}\label{Lem: minor equivalent}
For every class~$\calC$ of graphs we have $\Min(\calC) \equiv_\IN \calC$\,.
\end{lem}

\section{Strictness of the hierarchy}   
\label{Sect: hierarchy2}

In this section we prove that the hierarchy is strict.
Using the results of Section~\ref{Sect: tree decompositions and transductions}
we can characterise each level of the hierarchy
in terms of tree-width and its variants.
\begin{thm}\label{Thm: characterisation of pwd,twd}
Let $\calC \subseteq \STR[\Sigma]$.
\begin{enuma}
\item $\calC \sqsubseteq_\IN \calP$ \quad\iff\quad $\pwd(\calC) < \infty$.
\item $\calC \sqsubseteq_\IN \calT_\omega$ \quad\iff\quad $\twd(\calC) < \infty$.
\item $\calC \sqsubseteq_\IN \calT_n$ \quad\iff\quad $\twd_n(\calC) < \infty$.
\end{enuma}
\end{thm}
\begin{Proof}
In each case $(\Leftarrow)$ follows from Lemma~\ref{Lem: tree decomposable implies interpretable}
and $(\Rightarrow)$ follows from Theorem~\ref{Thm: interpretable implies tree decomposable}.
\end{Proof}
\begin{cor}
Let $\calC$~be a class of $\Sigma$-structures.
\begin{enuma}
\item $\pwd(\calC) = \infty$ implies $\calT_\omega \sqsubseteq_\IN \calC$.
\item $\twd(\calC) = \infty$ implies $\calG \sqsubseteq_\IN \calC$.
\end{enuma}
\end{cor}
\begin{Proof}
(a) Suppose that $\pwd(\calC) = \infty$.
Then Theorem~\ref{Thm: excluded tree theorem} implies
that $\calT_\omega \subseteq \Min(\GF(\calC))$.
Hence, the claim follows from
Lemmas \ref{Lem: Gaifman graph interpretable}~and~\ref{Lem: minor equivalent}.

(b) Suppose that $\twd(\calC) = \infty$.
Then Theorem~\ref{Thm: excluded grid theorem} implies
that $\calG \subseteq \Min(\GF(\calC))$.
As in~(a), the claim follows from
Lemmas \ref{Lem: Gaifman graph interpretable}~and~\ref{Lem: minor equivalent}.
\end{Proof}
\begin{cor}\label{Cor: nothing between P,T,G}
Let $\calC \subseteq \STR[\Sigma]$.
\begin{enuma}
\item $\calC \nsqsubseteq_\IN \calP$ implies $\calT_\omega \sqsubseteq_\IN \calC$.
\item $\calC \nsqsubseteq_\IN \calT_\omega$ implies $\calG \equiv_\IN \calC$.
\end{enuma}
\end{cor}

In particular, it follows that the upper part of the hierarchy is strict\?:
\begin{cor}\label{Cor: upper part}
$\calP \Sqsubset_\IN \calT_\omega \Sqsubset_\IN \calG$
\end{cor}
\begin{Proof}
Both assertions follow from
Theorem~\ref{Thm: characterisation of pwd,twd} and Corollary~\ref{Cor: nothing between P,T,G}.

For the first one, note that
we have $\calP \sqsubseteq_\IN \calT_\omega$ since $\calP \subseteq \Min(\calT_\omega)$.
Conversely, $\pwd(\calT_\omega) = \infty$ implies,
by Theorem~\ref{Thm: characterisation of pwd,twd}\,(a),
that $\calT_\omega \nsqsubseteq_\IN \calP$.
Hence, $\calP \sqsubset_\IN \calT_\omega$.
Finally, if $\calC \sqsubset_\IN \calT_\omega$ then $\calC \sqsubseteq_\IN \calP$,
by Corollary~\ref{Cor: nothing between P,T,G}\,(a).
Consequently, we have $\calP \Sqsubset_\IN \calT_\omega$.

Similarly, the fact that $\calT_\omega \sqsubseteq_\IN \calG$ follows
from Lemma~\ref{Lem: G is top}.
Since $\twd(\calG) = \infty$, Theorem~\ref{Thm: characterisation of pwd,twd}\,(b) implies
that $\calT_\omega \sqsubset_\IN \calG$.
Finally, we obtain $\calT_\omega \Sqsubset_\IN \calG$ by
Corollary~\ref{Cor: nothing between P,T,G}\,(b).
\end{Proof}

Let us turn to the lower part. We start with two technical lemmas.
\begin{defi}
Let $\frakT = \langle T,{\leq}\rangle$ be an order-tree.
Vertices $v_0,\dots,v_{m-1}$ of~$T$ are \emph{horizontally related via} a vertex~$w$
if all~$v_i$ are at the same level of the tree and $v_i \sqcap v_k = w$,
for all $0 \leq i < k < m$.
\end{defi}

\begin{lem}\label{Lem: large degree is preserved}
Let $\frakT$ be a coloured order-tree of height~$n$,
and suppose that $\tau$~is a parameterless $k$-copying $\MSO$-trans\-duc\-tion of rank~$r$
such that $\tau(\frakT)$ is a successor-tree of height at most $n+1$.
Consider vertices $v_0,\dots,v_{m-1}$ of~$\frakT$ that are horizontally related via~$w$
and fix some number $l < k$.
Let $x_i$~be the successor of~$w$ with $x_i \preceq v_i$.
If, for all $i,j < m$, we have
\begin{align*}
  \MTh_{r+2n+1}(\frakT_{x_i},v_i) &= \MTh_{r+2n+1}(\frakT_{x_j}, v_j)\,,
\end{align*}
then at least $m-1$ elements of the set
$\{\langle v_0,l\rangle,\dots,\langle v_{m-1},l\rangle\}$
\textup(these are elements of the domain of $\tau(\frakT)$\textup)
are horizontally related in $\tau(\frakT)$.
\end{lem}
\begin{Proof}
Let $\varphi_{ss'}(x,y)$ be the formula defining the successor relation in $\tau(\frakT)$
between vertices of the form $\langle x,s\rangle$ and $\langle y,s'\rangle$.
By assumption the rank of~$\varphi_{ss'}(x,y)$ is at most~$r$.

First, note that a vertex $\langle v,l\rangle$ is on level~$h$ in $\tau(\frakT)$
if and only if there are indices $s_0,\dots,s_{h-1} < k$ such that
\begin{align*}
  &\frakT \models \psi_{s_0\dots s_{h-1}}(v)
\intertext{where the formula}
  &\psi_{s_0\dots s_{h-1}}(v) :=
    \begin{aligned}[t]
      \exists y_0\cdots\exists y_{h-1}\Bigl[\Land_{i < h-1} \varphi_{s_is_{i+1}}(y_i,y_{i+1})
      &\land \varphi_{s_{h-1}l}(y_{h-1},v) \\
      &\land \neg \exists z\Lor_{s < k} \varphi_{ss_0}(z,y_0)\Bigr]
    \end{aligned}
\end{align*}
expresses that there exists a path of the form
$\langle y_0,s_0\rangle,\dots,\langle y_{h-1},s_{h-1}\rangle,\langle v,l\rangle$
from the root $\langle y_0,s_0\rangle$ of~$\tau(\frakT)$ to the vertex $\langle v,l\rangle$.
By assumption on~$v_i$ and Lemma~\ref{Lem: composition for trees}, we have
\begin{align*}
  \MTh_{r+2n+1}(\frakT,v_i) = \MTh_{r+2n+1}(\frakT,v_j)\,,
  \quad\text{for all } i,j\,.
\end{align*}
Since the rank of $\psi_{s_0\dots s_{h-1}}$ is $h+r+1 \leq r+2n+1$ it follows that
\begin{align*}
  \frakT \models \psi_{s_0\dots s_{h-1}}(v_i)
  \quad\iff\quad
  \frakT \models \psi_{s_0\dots s_{h-1}}(v_j)\,.
\end{align*}
Hence, all vertices $\langle v_0,l\rangle,\dots,\langle v_{m-1},l\rangle$ are on the same level~$h$
in $\tau(\frakT)$.
We prove by induction on~$h$ that
\begin{align*}
\prefixtext{$(*)$}
  \MTh_{r+n+h+1}(\frakT_{x_i},v_i) = \MTh_{r+n+h+1}(\frakT_{x_j}, v_j)
\end{align*}
implies that all but at most one of $\langle v_0,l\rangle,\dots,\langle v_{m-1},l\rangle$
are horizontally related.

Let $\langle u_i,s_i\rangle$~be the predecessor of $\langle v_i,l\rangle$ in $\tau(\frakT)$,
that is,
\begin{align*}
  \frakT \models \varphi_{s_il}(u_i,v_i)\,.
\end{align*}
We distinguish two cases.

First suppose that $u_0 \sqcap v_0 \preceq w$ in~$\frakT$.
By~$(*)$ and Lemma~\ref{Lem: composition for trees}, we have
\begin{align*}
  \MTh_{r+n+h+1}(\frakT,u_0,v_0) = \MTh_{r+n+h+1}(\frakT,u_0,v_i)\,,
\end{align*}
for all~$i$ such that $u_0 \sqcap v_i \preceq w$.
Note that there can be at most one index~$i$ that does not satisfy this condition
since, if we had $u_0 \sqcap v_i \succeq w$ and $u_0 \sqcap v_j \succeq w$, for $i \neq j$,
then we would have $v_i \sqcap v_j \succ w$ and $v_0,\dots,v_{m-1}$ would not be
horizontally related via~$w$.
It follows that
\begin{align*}
  \frakT \models \varphi_{s_0l}(u_0,v_0)
  \qtextq{implies}
  \frakT \models \varphi_{s_0l}(u_0,v_i)\,,
  \quad\text{for all } i \text{ as above}\,.
\end{align*}
Hence, $\langle u_0,s_0\rangle$ is the common predecessor of all the $\langle v_i,l\rangle$,
except for possibly one of them.
(For an index~$i$ with $u_0 \sqcap v_i \succeq w$ our composition
argument does not work since in that case $(*)$~does not imply that the theories of
$(\frakT,u_0,v_0)$ and $(\frakT,u_0,v_i)$ coincide.)

It remains to consider the case that $w \prec u_0 \sqcap v_0$.
Setting
\begin{align*}
  \eta_{u_i} := \Land \MTh_{r+n+h-1+1}(\frakT_{x_i},u_i)
\end{align*}
we have
\begin{align*}
  \frakT_{x_0} \models \exists z[\abs{z} = \abs{u_0} \land \eta_{u_0}(z) \land \varphi_{s_0l}(z,v_0)]\,.
\end{align*}
Since the rank of this formula is $r+n+h+1$ it follows that
\begin{align*}
  \frakT_{x_i} \models \exists z[\abs{z} = \abs{u_0} \land \eta_{u_0}(z) \land \varphi_{s_0l}(z,v_i)]\,,
  \quad\text{for all } i < m\,.
\end{align*}
Consequently, we have $\abs{u_i} = \abs{u_0}$, for all~$i$,
and $u_0,\dots,u_{m-1}$ are horizontally related via~$w$.
Since the vertices $\langle u_0,s_0\rangle,\dots,\langle u_{m-1},s_0\rangle$
are on level $h-1$ in $\tau(\frakT)$,
we can apply the induction hypothesis and it follows that all but at most one of
then are horizontally related via some vertex~$w'$.
Therefore, the same holds for their
successors $\langle v_0,l\rangle,\dots,\langle v_{m-1},l\rangle$.
\end{Proof}

\begin{defi}
We denote by $B(n,k,c)$ the number of
functions $m^{<n} \to \PSet([c])$ with $m \leq k$.
Intuitively, each such function corresponds to
a vertex-colouring of the tree $m^{<n}$ with $c$~colours.
\end{defi}
\begin{lem}
For $n \geq 1$ and $k \geq 2$, we have
\begin{align*}
  2^{ck^{n-1}} \leq B(n,k,c) \leq k2^{2ck^{n-1}}.
\end{align*}
\end{lem}
\begin{Proof}
For $m \geq 2$, we have
\begin{align*}
  m^{n-1} \leq m^{n-1} + \sum_{i < n-1} m^i = m^{n-1} + \frac{m^{n-1}-1}{m-1} \leq 2m^{n-1}.
\end{align*}
Since $\bigabs{[m]^{<n}} = \sum_{i < n} m^i = m^{n-1} + \sum_{i < n-1} m^i$ it follows that
\begin{align*}
  m^{n-1} \leq \bigabs{[m]^{<n}} \leq 2m^{n-1}.
\end{align*}
Therefore, we can bound
\begin{align*}
  B(n,k,c) = 2^{cn} + \sum_{m = 2}^k 2^{c\abs{[m]^{<n}}}
\end{align*}
from above by
\begin{align*}
  B(n,k,c) \leq 2^{cn} + \sum_{m = 2}^k 2^{c2m^{n-1}} \leq k2^{2ck^{n-1}}
\end{align*}
and from below by
\begin{align*}
  B(n,k,c) \geq 2^{cn} + \sum_{m = 2}^k 2^{cm^{n-1}} \geq 2^{ck^{n-1}}.
\end{align*}
\upqed
\end{Proof}

\begin{thm}\label{Lem: lower part strict}
$\calT_n \sqsubset_\IN \calT_{n+1}$
\end{thm}
\begin{Proof}
For a contradiction, suppose that there exists a transduction~$\tau$ such that
$(\calT_{n+1})_\IN \subseteq \tau((\calT_n)_\IN)$.
Let $\calT_n^\ord$~be the class of all order-trees corresponding to successor-trees in~$\calT_n$,
and let $\calT_{n+1}^{\mathrm{col}} := \mathrm{exp}_1(\calT_{n+1})$
be the class of all coloured successor-trees with one colour
whose underlying tree is in~$\calT_{n+1}$.
Since the suc\-cessor-trees in~$\calT_n$ are $1$-sparse we can construct an
$\MSO$-trans\-duc\-tion~$\sigma_0$ such that $(\calT_n)_\IN \subseteq \sigma_0(\calT_n)$.
Since $\calT_n^\ord \equiv \calT_n$ we can combine $\tau$~and~$\sigma_0$
to a transduction~$\sigma$ such that $\calT_{n+1}^{\mathrm{col}} \subseteq \sigma(\calT_n^\ord)$.
By Lemma~\ref{Lem: large degree is preserved}, it follows that there is some
constant~$d$ such that every tree $\frakT \in \sigma(\calT_n^\ord)$ with out-degree at most~$k$
is of the form $\sigma(\frakT')$ where $\frakT' \in \calT_n^\ord$ has out-degree at most~$dk$.
(The out-degree of an order-tree is the out-degree of the corresponding successor-tree.)
Suppose that $\sigma$~uses $c$~parameters. There are
\begin{align*}
  B(n,dk,c) \leq dk2^{2c(dk)^{n-1}}
\end{align*}
colourings of trees in $\calT_n^\ord$ with out-degree at most~$dk$.
On the other hand, there are
\begin{align*}
  B(n+1,k,1) \geq 2^{k^n}
\end{align*}
trees in $\calT_{n+1}^{\mathrm{col}}$ with out-degree at most~$k$.
For large~$k$ it follows that
\begin{align*}
  B(n,dk,c) \leq dk2^{2cd^{n-1}k^{n-1}}
  < 2^{k^n} = B(n+1,k,1)\,.
\end{align*}
Consequently, there is some tree in $\calT_{n+1}^{\mathrm{col}}$ that is not the image
of a tree in~$\calT_n^\ord$. A~contradiction.
\end{Proof}

\section{Completeness of the hierarchy}
\label{Sect: hierarchy3}

We have shown that the hierarchy presented in Theorem~\ref{Thm: hierarchy}
is strict. To conclude the proof of the theorem
it therefore remains to show that there are no additional classes.
We have already seen in Corollary~\ref{Cor: upper part}
that $\calT_\omega$~and~$\calG$
are the only classes above $\calP$.
Next we shall prove that $\calT_n \Sqsubset \calT_{n+1}$.
\begin{lem}\label{Lem: embeddable trees means interpretable}
Let $\calC$~be a class of structures. If, for every number $m \in \bbN$,
there exists a structure $\frakA \in \calC$ such that
we can embed $m^{<n+1}$ into every tree underlying a tree decomposition
$(U_v)_{v \in T}$ of~$\frakA$ of width~$k$,
then $\calT_{n+1} \sqsubseteq_\IN \calC$.
\end{lem}
\begin{Proof}
By Lemma~\ref{Lem: existence of strict tree decompositions}, it follows that,
for every $m \in \bbN$, there is a structure in~$\calC$ with a strict tree decomposition
of width at most~$k$ and with an underlying tree~$T$ into which we can embed the tree~$m^{<n+1}$.
According to Proposition~\ref{Prop: strict tree decompositions definable}
there is a transduction mapping~$\calC$ to
the class of trees underlying these strict tree decompositions.
Hence, there exists a class $\calK \sqsubseteq_\IN \calC$ of successor-trees
containing, for every $m \in \bbN$, some tree into which we can embed $m^{<n+1}$.
Hence, $\calT_{n+1} \subseteq \Min(\calK) \sqsubseteq_\IN \calK \sqsubseteq_\IN \calC$.
\end{Proof}

\begin{thm}\label{Thm: covering for lower part}
Let $\calC$~be a class of structures.
Then $\calT_{n+1} \nsqsubseteq_\IN \calC$ implies $\twd_n(\calC) < \infty$.
\end{thm}
\begin{Proof}
Suppose that $\calT_{n+1} \nsqsubseteq_\IN \calC$.
Then $\calP \nsqsubseteq_\IN \calC$ since $\calT_{n+1} \sqsubseteq_\IN \calP$.
By Lemmas \ref{Lem: Gaifman graph interpretable}~and~\ref{Lem: minor equivalent},
this implies that $\calP \nsubseteq \Min(\GF(\calC))$.
Therefore, we can find a path that is not in $\Min(\GF(\calC))$.
By Theorem~\ref{Thm: excluded path theorem}, it follows that
there are numbers $k,l \in \bbN$ such that
$\twd_l(\calC) < k$.

By induction on~$l$, we prove that
$\twd_l(\calC) < \infty$ implies $\twd_n(\calC) < \infty$.
For $l \leq n$, there is nothing to do.
For $l > n$, we have
$\calT_{n+1} \sqsubseteq_\IN \calT_l$, which implies that
$\calT_l \nsqsubseteq_\IN \calC$.
Consequently, it follows by Lemma~\ref{Lem: embeddable trees means interpretable} and
Corollary~\ref{Cor: bound on twd-n}
that $\twd_{l-1}(\calC) < \infty$.
By induction hypothesis, the result follows.
\end{Proof}
By Lemma~\ref{Lem: tree decomposable implies interpretable}
we obtain the following results.
\begin{cor}\label{Cor: nothing between Tn and Tn+1}
If $\calT_{n+1} \nsqsubseteq_\IN \calC$ then $\calC \sqsubseteq_\IN \calT_n$.
\end{cor}
\begin{cor}\label{Cor: lower part covers}
$\calT_n \Sqsubset_\IN \calT_{n+1}$.
\end{cor}

To conclude the proof of Theorem~\ref{Thm: hierarchy} it remains to show
that there are no classes between the lower part of the hierarchy and its upper part.
\begin{lem}\label{Lem: calT_n dense below calP}
Let $\calC$~be a class of structures.
If $\calT_n \sqsubseteq_\IN \calC$, for all $n \in \bbN$, then $\calP \sqsubseteq_\IN \calC$.
\end{lem}
\begin{Proof}
We show the contraposition.
Suppose that $\calP \nsqsubseteq_\IN \calC$.
We have to show that there is some~$n$
such that $\calT_n \nsqsubseteq_\IN \calC$.
As in the proof of Theorem~\ref{Thm: covering for lower part}
it follows that there are numbers $k,l \in \bbN$ such that $\twd_l(\calC) < k$.
Hence, we can use Lemma~\ref{Lem: tree decomposable implies interpretable}
to obtain a transduction~$\tau$ witnessing that $\calC \sqsubseteq_\IN \calT_l$.
By Theorem~\ref{Lem: lower part strict}
we have $\calT_{l+1} \nsqsubseteq_\IN \calT_l$.
It follows that $\calT_{l+1} \nsqsubseteq_\IN \calC$, as desired.
\end{Proof}

\begin{cor}\label{Cor: noting between lower and upper part}
If $\calC \sqsubset_\IN \calP$ then there is some $n \in \bbN$ such that $\calC \sqsubseteq_\IN \calT_n$.
\end{cor}
\begin{Proof}
By Lemma~\ref{Lem: calT_n dense below calP}, there is some $n \in \bbN$ such that
$\calT_{n+1} \nsqsubseteq_\IN \calC$. Hence, Corollary~\ref{Cor: nothing between Tn and Tn+1}
implies that $\calC \sqsubseteq_\IN \calT_n$.
\end{Proof}

Together, Corollaries \ref{Cor: upper part},~\ref{Cor: lower part covers},
and~\ref{Cor: noting between lower and upper part}
(and the fact that every class~$\calC$ satisfies $\emptyset \sqsubseteq_\IN \calC \sqsubseteq_\IN \calG$)
show that every class of $\Sigma$-structures
is $\equiv_\IN$-equivalent to some of the classes
in Theorem~\ref{Thm: hierarchy}.
This completes the proof of this theorem.

\section{Prospects and conclusion}   
\label{Sect: conclusion}

Above we have obtained a complete description of the transduction hierarchy for classes
of finite incidence structures.
The most surprising result is that the hierarchy is linear.
At this point there are at least three natural directions
in which to proceed.
\begin{enumi}
\item We can study the hierarchy for classes of structures,
  instead of classes of incidence structures.
\item We can consider the hierarchy for classes of infinite structures.
\item We can replace $\MSO$ by a different logic.
\end{enumi}

An answer to~(ii) seems within reach, at least if we restrict our attention to countable
structures. Although the resulting hierarchy is no longer linear
we can adapt most of our techniques to this setting.
(For an example of nonlinearity, note that the class of all countable trees and
the class of all finite grids are incomparable.)

Concerning question~(iii), let us remark that all results above go through if we use $\CMSO$
instead of~$\MSO$. We only need the right definition of rank for $\CMSO$.
In the proof of Theorem~\ref{Lem: lower part strict}
we needed the fact that there are only finitely many theories of bounded rank.
We can ensure this for $\CMSO$ by defining the rank as the least number~$n$
such that
\begin{itemize}
\item the nesting depth of quantifiers is at most~$n$ and
\item in every cardinality predicate $\abs{X} \equiv k \pmod m$ we have $m \leq n$.
\end{itemize}
One can check that, with this definition of rank, the proof
of Theorem~\ref{Thm: composition for unions} also goes through for $\CMSO$.

For logics much weaker than $\MSO$, on the other hand, it seems unrealistic
to hope for a complete description of the corresponding transduction hierarchy.
For instance, a related hierarchy for first-order logic was investigated
by Mycielski, Pudl\'ak, and Stern in~\cite{MycielskiPudlakStern90}.
The results they obtain indicate that the structure of the resulting hierarchy
is very complicated.

Finally, let us address question~(i).
When using transductions between structures instead of their incidence structures,
we can transfer some of the above results to the corresponding hierarchy.
But we presently have no complete description since we miss some of the corresponding excluded minor results.

\begin{lem}\label{Lem: relation between ad and in}
Let $\calC,\calS \subseteq \STR$ and suppose that $\calS$~is $k$-sparse.
\begin{enuma}
\item $\calC \sqsubseteq_\IN \calS$ implies $\calC \sqsubseteq \calS$.
\item $\calS \sqsubseteq \calC$ implies $\calS \sqsubseteq_\IN \calC$.
\end{enuma}
\end{lem}
\begin{Proof}
There is a transduction~$\varrho$ such that $\calC = \varrho(\calC_\IN)$.
Since $\calS$~is $k$-sparse we can also find a transduction~$\sigma$ such that
$\calS_\IN = \sigma(\calS)$.
Consequently,
\begin{alignat*}{-1}
  \calC_\IN &\subseteq \tau(\calS_\IN) &&\qtextq{implies}
  &\calC &\subseteq (\varrho \circ \tau \circ \sigma)(\calS)\,, \\
\prefixtext{and}
  \calS &\subseteq \tau(\calC) &&\qtextq{implies}
  &\calS_\IN &\subseteq (\sigma \circ \tau \circ \varrho)(\calC_\IN)\,.
\end{alignat*}
\upqed
\end{Proof}

\begin{thm}\label{Thm: hierarchy for MSO}
We have the following hierarchy\?:
\begin{align*}
  \emptyset \sqsubset \calT_0 \sqsubset \calT_1
  \sqsubset \dots \sqsubset \calT_n \dots \sqsubset \calP
  \sqsubset \calT_\omega
  \sqsubset \calG \equiv \STR[\Sigma]
\end{align*}
\end{thm}
\begin{Proof}
Note that all classes in Theorem~\ref{Thm: hierarchy for MSO} are $2$-sparse. For $2$-sparse
classes $\calC$ and~$\calK$, Lemma~\ref{Lem: relation between ad and in} implies that
\begin{align*}
  \calC \sqsubseteq \calK \quad\iff\quad
  \calC \sqsubseteq_\IN \calK\,.
\end{align*}
Consequently, the result follows from Theorem~\ref{Thm: hierarchy}.
\end{Proof}

\begin{oprob}
Is there any class $\calC \subseteq \STR[\Sigma]$ which is not $\equiv$-equivalent
to some class in the above hierarchy\??
\end{oprob}

\begin{rem}
If we only consider classes of graphs
and if we use $\CMSO$-transductions instead of $\MSO$-transductions,
then the following result
can be used as replacement of
Theorem~\ref{Thm: excluded grid theorem}\?:
\end{rem}
\begin{thm}[\cite{CourcelleOum04}]
Let $\calC$~be a class of graphs with unbounded clique-width.
There exists a $\CMSO$-transduction~$\tau$ with
$\calG \subseteq \tau(\calC)$.
\end{thm}
\noindent
This eliminates some possibilities for intermediate classes of graphs in the hierarchy
of Theorem~\ref{Thm: hierarchy for MSO},
but to complete the picture we still need analogues of
Proposition~\ref{Prop: strict tree decompositions definable}
and of Theorems \ref{Thm: excluded tree theorem}~and~\ref{Thm: excluded path theorem}.
Furthermore, the techniques of~\cite{CourcelleOum04} are specific to graphs
(or, more generally, to relational structures where all relations are binary).
Even with the results of~\cite{CourcelleOum04} one cannot exclude the existence of a class~$\calC$
of arbitrary relational structures strictly between
$\calT_\omega$~and~$\calG$ in the $\CMSO$-transduction hierarchy.

Let us make a final comment about relational structures.
An incidence structure~$\frakA_\IN$ can be seen as a bipartite labelled directed graph
(see the remark after Definition~\ref{Def: incidence structure}).
Furthermore, it is $1$-sparse.
Hence, our results use tools from graph theory,
in particular those of \cite{RobertsonSeymour83,RobertsonSeymour86,CourcelleOum04}.
However, there is currently no encoding of relational structures as labelled graphs
that could help to solve question~(i) above.

\end{document}